\documentclass{amsproc}
 \usepackage{graphicx}
\input amssymb.sty
\def \\ { \cr }
\def \R{{\bf R}}
\def \N{{\bf N}}

\def \P{{\bf P}}
\def \C{{\bf C}}
\def \Q{{\bf Q}}

\def \i{{\infty}}

\def \th{{\theta}}

\def\z{{\zeta}}
\def\l{{\lambda}}
\newcommand {\cc}{\check}

\newtheorem{definition}{Definition}[section]

\newtheorem{theorem}[definition]{Theorem}

\newtheorem{proposition}[definition]{Proposition}
\newtheorem{lemma}[definition]{Lemma}

\def \og{\leavevmode\raise.3ex\hbox{$\scriptscriptstyle\langle\!\langle $}}
\def \fg{\leavevmode\raise.3ex\hbox{$\scriptscriptstyle\rangle\!\rangle \ $}}

\begin{document}

\title{Flutuations of L\'evy processes and Scattering Theory}
\author{Sonia Fourati}

\address{I.N.S.A. de Rouen. Place Emile Blondel. 76130 Mont Saint-Aignan, France}
\address{L.P.M.A. Universit\'e Paris VI. Case courrier 188. 4 place Jussieu 75252 Paris Cedex 05}
\email{fourati@ccr.jussieu.fr}
\date{ 21 December 2006}

\begin{abstract} We establish a connection between  the inverse scattering problem and the determination of the distribution of the position of a L\'evy process at the exit time of a bounded interval in term of its L\'evy exponent.

\end{abstract}

\maketitle

\section {Introduction}

It is well known that the fluctuation properties of a L\'evy process are intimately related to the Wiener-Hopf factorization of its L\'evy exponent $\phi$. The observation (initiated by Spitzer \cite{Sp} for random walks) that  the factorization of $(q+\phi(iu))^{-1}$, where $q$ is a positive constant, can be interpreted as the 
independence of the past and pre minimum  parts of the L\'evy process killed at an independent exponential time. The Wiener-Hopf factors yield the  distributions of the minimum and the maximum of the killed process. Also the  problem of exit from a semi-infinite interval which consists in determining the  joint distribution of the exit time of an interval $]-\infty,x]$ and of the position at the exit time,
 can be reduced to finding  the Wiener-Hopf factorization, see \cite{B}\cite {D}\cite {S}.
 In the theory of analytic functions, Wiener-Hopf factorization is the simplest of a large class of factorization
problems known as Riemann-Hilbert problems. 
In this paper we establish a connection between the  joint distribution of the maximum, minimum and final value of the L\'evy process killed at an independent exponential time and a certain Riemann-Hilbert factorization problem. The determination of this joint distribution allows to compute the joint law of the position and the time of the exit of a L\'evy process from a bounded interval.
  More precisely, we prove that this  problem reduces to  the  factorization of the matrix $$\left( \begin{array} {clcr} 0 & -1\\1 & \phi(iu ) +q\end{array} \right)$$
 where $\phi$ is the L\'evy exponent of the process, into a product 
 $A(x,iu)B(x, iu)$ where  $A$ and $B$ are matrices which are boundary values, on the imaginary axis, of analytic functions, defined  respectively  on the left and on the  right complex half planes, satisfying a normalization at  infinity involving a positive real parameter $x$. In order to  obtain this result we establish that a set of functions, defined from
 the Laplace transforms of random variables involving the maximum and minimum processes associated with the L\'evy process, 
 satisfy a certain system of integral equations. We show that this system of integral equation is analogous to the basic differential system appearing in   scattering theory on the real  line.
 Our system does not reduce to the usual problem of scattering theory, since the associated potential is very singular compared to the potentials considered in this theory. The analogy however is sufficiently good that we can apply similar arguments as in the work of Shabat \cite{Sh}, and show the equivalence of our system of integral equations with a  Riemann-Hilbert problem.
 
 This paper is organized as follows. In section 2 we recall some basic facts on L\'evy processes and the Wiener-Hopf factorization. Then we introduce the main functions of our work and state the two main theorems. The first one gives a system of integral equations satisfied by the functions, and the second one states that this system is equivalent to a certain Riemann-Hilbert problem. We explain why the first theorem is related to the differential equation and the direct problem of the scattering theory while the second one is 
 related to the inverse scattering problem. In section 3 we explain the connection between the main results and the exit problem from an interval and related questions. Section 4 contains a preliminary result on the conditional independence of the pre and post minimum process, knowing the amplitude. In section 5 we establish the first equation of the integral system. Section 6 deals with the second group of equations. In this part we introduce two Markov chains, built from the successive minima and maxima of the process, which play a key role in the proof. In section 7 we establish the second theorem, about the Riemann-Hilbert problem, by adapting the arguments of Shabat \cite{Sh}
 to our setting.  
 We give some probabilistic interpretations of  
  the factorization in the Riemann-Hilbert problem in section 8 in terms of Wiener-Hopf factorization of certain auxiliary L\'evy processes.   
   Finally, section 9 is devoted to apply our results to stable processes and to L\'evy processes without positive jumps. This allows us to give some precisions on results obtained by Rogozin \cite {R} for the stable case and by Takacs \cite {T} for  L\'evy processes without positive jumps.

\section{Notations and main results} 

Let $(\Omega,{\mathcal F})$ be  the  space of functions defined on $[0,+\i[$ with 
 values in $\R\cup \{\delta\}$ where $\delta$ is a cemetery point, and let $X$ denote the canonical process $X_t(\omega)=\omega(t)$.
  In this  paper $\P$ will be  the law 
 on   $(\Omega,{\mathcal F})$ of a L\'evy process started at $0$ with L\'evy exponent $\phi$. More precisely, we have :
$$ \P(e^{-iu X_t})= e^{-t \phi(iu)}\qquad (iu\in i\R)$$

\subsection{Some facts on Wiener-Hopf factorization}
We start by recalling some standard facts on  Wiener-Hopf factorization and fluctuations of L\'evy processes for which we refer to [B] chapter 6.
Let  $S_t$ et $I_t$ be the past maximum and past minimum processes, namely :

$$S_t:=\sup \{X_s, 0\leq s\leq t\}\qquad I_t :=\inf  \{X_s, 0\leq s\leq t\}$$

We introduce now local times at $0$  of the reflected processes  $S-X$ et $X-I$ and the associated Wiener-Hopf factors. The definition of these local times 
depends on the regularity of $[0,+\i[$ or $]-\i,0]$ for the L\'evy process.

If  $]0,+\i[$ is regular, i.e. $T^0=\inf\{t>0, X_t> 0\}=0$  a.s.  (resp. $]-\i,0[$ is regular $T_0=\inf\{t>0, X_t<0\}=0$ a.s.) then 
$0$ is a regular point for the Markov process $X-S$ (resp. $X-I$) and $L^s$ (resp. $L^i$) denotes any local time at $0$ of this process.
In this case, $t\mapsto L^s[0,t]$ (resp. $t\mapsto L^i[0,t]$) is an increasing continuous process, we denote by $L^{s,-1}$ (resp.$L^{i,-1}$) its right continuous inverse. The  pair $(L^{s,-1}, S_{L^{s,-1}})$ (resp. $(L^{i,-1},I_{L^{i,-1}})$) is a bi-variate L\'evy process which may have a finite life time
 if $\lim_{t\to +\i}  X_t=-\i$ a.s. (resp. $\lim_{t\to +\i}  X_t=+\i$ a.s.).
The Wiener-Hopf factors are the L\'evy exponents of this process, more precisely

$$\P(e^{-\l S_{L^{s,-1}_t}-  qL^{s,-1}_t}; t<L^s_{\i})=:e^{-t\psi_q(\l)}\qquad \Re(\l)\geq 0 ,q\in [0,+\i[$$
Respectively,
$$\P(e^{-\l I_{L^{i,-1}_t}-  qL^{i,-1}_t}; t<L^i_{\i})=:e^{-t\cc \psi_q(\l)}\qquad \Re(\l)\leq 0,q\in [0,+\i[$$

If  $[0,+\i[$ is irregular (this means that the time $T^{0-}:=\inf\{t>0, X_t\geq 0\}$ is  positive a.s.) , (respectively if $]-\i,0]$ is irregular , $T_{0-}:=\inf\{t>0, X_t\leq 0\}$ is  positive a.s.) then  the set  $\{t; S_t=X_t\}$ (resp. $\{t; X_t=I_t\}$) is a.s. discrete   and  the local time $L^s$ (respectively, $L^i$) is the random point measure $$L^s(dt):=\sum_u1_{X_u=S_u}\delta_u(dt)$$
Respectively, $$L^i(dt):=\sum_u1_{X_u=I_u}\delta_u(dt)$$

The Wiener-Hopf factors are  
$$\psi_q(\l):=1-\P(e^{-\l X_{T^{0-}}-qT^{0-}}; T^{0-}<+\i) \qquad \Re(\l)\geq 0, q \in [0,+\i[ $$
Respectively,
$$\cc \psi_q(\l):=1-\P(e^{-\l X_{T_{0-}}-qT_{0-}}; T_{0-}<+\i) \qquad \Re(\l)\leq 0, q \in [0,+\i[ $$
Notice that in this case $\psi$ (resp. $\cc \psi$) is the L\'evy exponent of a compound Poisson process and that 
according to proposition 4 of [B] chapter 6, $X_{T^{0-}}>0$ a.s.  (resp. $X_{T_{0-}}<0$ a.s. ), thus the times $T^{0-}$ and $T^0:=\inf \{t, X_t>0\}$ are equal a.s. (resp.
$T_{0-}=T_0:=\inf\{t; X_t<0\}$ a.s.).

Finaly, if neither condition is fulfilled, then $\P$ is   the law a compound Poisson process. In this case, for reasons which will appear later, it is 
necessary to use a dissymetric definition of local times : We denote by 
$L^s$ the random measure 
$$L^s(dt):=1_{S_t=X_t}dt$$
And $L^i$ will be the random point measure :
$$L^i(dt):=\delta_0(dt)+\sum_{u>0}1_{I_{u^-}>I_u}1_{X_u=I_u}\delta_u(dt)$$
As above the Wiener-Hopf factor $\psi_q(\l)$ is the L\'evy exponent of the   bi-variate L\'evy process  
$(L^{s,-1}, S_{L^{s,-1}})$, namely :
$$\P(e^{-\l S_{L^{s,-1}_t}-  qL^{s,-1}_t})=:e^{-t\psi_q(\l)}\qquad \Re(\l)\geq 0,q\in [0,+\i[$$
The Wiener-Hopf factor $\cc \psi_q(\l)$ is the fonction :
$$\cc\psi_q(\l):=1-\P(e^{-\l X_{T_0}-qT_0}; T_0<+\i)\qquad T_0=\inf\{t, X_t<0\}$$

Note that in all  cases, one has :

 $${1\over \psi_q(\l)}=\P(\int_{[0,+\i[}e^{-\l S_t-qt}L^s(dt))\qquad \Re(\l)\geq 0, q\in [0,+\i[, q\Re(\l)\not=0$$
 and 
 $${1\over \cc\psi_q(\l)}=\P(\int_{[0,+\i[}e^{-\l I_t-qt}L^i(dt))\qquad\Re(\l)\leq 0, q\in [0,+\i[, q\Re(\l)\not=0$$
It is  possible to normalize the local times so that the following Wiener-Hopf factorization holds
  (see e.g. \cite{B},\cite{S}), and we shall assume in the sequel that it is the case. For the compound Poisson process this follows from the convention we have choosen. 
\begin{proposition}\label{WH}
The pair  $(\psi_q(\l),\cc \psi_q(\l))$ satisfies the following Wiener-Hopf  identity 
$$\cc\psi_q(iu)\psi_q(iu)=\phi(iu)+q\qquad iu\in i\R, q\in [0,+\i[$$
\end{proposition}
\subsection{The main functions}
In next proposition we define  the so called excursions measures $N$ and $\cc N$ associated to local times $L^i$ and $L^s$ by the   compensation formula (see for exemple chapter 4 of [B]). 

\begin{proposition}\label {CF} {\bf Compensation formula} 

There exists a unique measure on $(\Omega,{\mathcal F})$, $N$ (resp. $\cc N$) such that 
$$\P\left(\sum_{]g,d[\in C}1_{(I_g-X_{(g-s)^-})_{s\geq 0}\in dw_1}1_{(X_{g+t}-I_g) _{0\leq t<d-g} \in dw_2}\right)
 =\P\left(\int_{[0,+\i[} 1_{(I_t-X_{(t-s)^-})_{ s\geq 0}\in dw_1}L^i(dt)\right) N( dw_2)$$
respectively, 
$$\P\left(\sum_{]g,d[ \in \cc C}1_{(S_g-X_{(g-s)^-})_{s\geq 0}\in dw_1} 1_{(X_{g+t}-S_g) _{0\leq t<d-g} \in dw_2}\right)=\P\left(\int_{[0,+\i[} 1_{(S_t-X_{(t-s)^-})_{s\geq 0}\in dw_1}L^s(dt)\right)\cc N( dw_2)$$
where $C$ (resp. $\cc C$)  is the set of connected components of the complement of  the  support of $L^i(dt)$ (resp. $L^s(dt)$).
\end{proposition}
Notice that if $ [0,+\i[$  (resp. $]-\i,0[$) is regular, then the state $0$ is a regular point of the Markov process $X-S$(resp. $X-I$) and $\cc N$ (resp. $ N$) is the usual excursion measure from $0$ of this process. If $[0,+\i[$  (resp. $]-\i, 0[$)  is irregular, then $\cc N$ (resp. $N$) is the distribution under $\P$ of the canonical process $X$ killed at time $T^0$ (resp. $T_0$).

We can now introduce the main functions of this paper. First define the random stopping times for every $x\in ]0,+\i[$ :

$$T^x:=\inf \{t; X_t>x\} \qquad T_x:=\inf \{t; X_t<-x\}$$
$$T^s _x:=\inf \{t ;X_t-S_t<-x\}\qquad T^x_i:=\inf \{t ;X_t-I_t>x\}$$

Define the right continuous left limited functions of $x$ ($x\in ]0,+\i[$) :

$$A_q(x,\l):=\P(\int_{[0,+\i[} 1_{S_t-I_t\leq x}e^{-\l S_t-qt}L^s(dt))\qquad  \l \in \C,  q \in [0,+\i[$$
$$\cc A_q(x,\l):=\P(\int_{[0,+\i[} 1_{S_t-I_t\leq x}e^{-\l I_t-qt}L^i(dt))\qquad  \l \in \C,  q\in [0,+\i[$$
$$C_q(x, \l):=N(e^{-\l X_{T^x}-qT^x}; T^x<+\i) \qquad   \Re(\l)\geq 0,  q\in [0,+\i[$$
$$\cc C_q(x,\l):=\cc N(e^{-\l X_{T_x}-qT_x}; T_x<+\i)  \qquad \Re(\l)\leq 0 , q\in [0,+\i[$$

If $[0,+\i[$ (resp. $ ]-\i,0[$) is regular, it is easy to check that the process  $(( L^{s,-1}_t; S_{L^{s,-1}_t}); 0\leq t<L^s_{T^s_x}$) (resp. $(( L^{i,-1}_t; I_{L^{i,-1}_t}); 0\leq t<L^i_{T^x_i})$)  is a killed L\'evy process and we denote by $B_q(x,\l)$ (resp. $\cc B_q(x,\l)$) its L\'evy exponent, more precisely :
$$e^{-t B_q(x,\l)}:=\P(e^{-\l S_{L^{s,-1}_t}-q L^{s,-1}_t }; t<L^s_{T^s_x})\qquad \Re(\l)\geq 0, q\in [0,+\i[$$
$$e^{-t \cc B_q(x,\l)}:=\P(e^{-\l I_{L^{i,-1}_t} -q L^{i,-1}_t}; t<L^i_{T^x_i})\qquad  \Re(\l)\leq 0, q\in [0,+\i[$$

If  $[0,+\i[$ (resp. $]-\i,0[$) is irregular, then
$$B_q(x,\l):=1- \P(e^{-\l S_{T^0}-q T^0};T^0<T^s_x)\qquad   \Re(\l)\geq 0, q\in [0,+\i[$$
respectively, 
$$\cc B_q(x,\l):=1- \P(e^{-\l I_{T_0}-q T_0}; T_0<T_i^x)\qquad  \Re(\l)\leq 0, q\in [0,+\i[\quad $$
Notice that $T^0<T^s_x$ if and only if $T^0<T_x$ (resp. $T_0<T_i^x$ if and only if $T_0<T^x$).  Notice also that  the function $(q,\l)\mapsto B_q(x,\l)$(resp. $(q,\l)\mapsto \cc B_q(x,\l)$) is the L\'evy exponent of a compound Poisson process.

In all cases one gets 
$${1\over B_q(x,\l)}=\P(\int_{[0,T^s_x[} e^{-\l S_t-qt} L^s(dt))\qquad  \Re(\l)\geq 0 ,q\in [0,+\i[$$
$${1\over \cc B_q(x,\l)}=\P(\int_{[0,T^x_i[} e^{-\l I_t-qt} L^i(dt))\qquad \Re(\l)\leq 0, q\in [0,+\i[$$

Using the previous functions, we now define the following  ones :
 
$$H_q(x):=A_q(x,0)=\P(\int_{[0,+\i[} 1_{S_t-I_t\leq x}e^{-qt}L^s(dt))$$
$$\cc H_q(x):=\cc A_q(x,0)=\P(\int_{[0,+\i[} 1_{S_t-I_t\leq x}e^{-qt}L^i(dt))$$
We shall denote $H_q(dx)$ and $\cc H_q (dx)$ the Stieltjes measures associated to  these increasing  functions.
Since $L^s[0,\varepsilon]$ and $L^i[0,\varepsilon]$ are positive for every $\varepsilon>0$, 
$H_q(x)$ and  $\cc H_q(x)$ do not vanish. Furthermore, one has :
$$H_q(x)\leq \P(\int_{[0,+\i[} 1_{S_t\leq x} e^{-qt} L^s(dt))\quad\hbox{and}\quad \cc H_q(x)\leq \P(\int_{[0,+\i[} 1_{-I_t\leq x} e^{-qt} L^i(dt))$$
For $q=0$, the right hand sides of the first  (resp. second) inequality is the so called renewal function of the subordinator with L\'evy exponent $\psi_0$ (resp. $\cc \psi_0$) (see [B] chapter 3). Therefore it is finite and 
$H_q(x)$ and $\cc H_q(x)$ are finite too. It is also true obviously when  $[0,+\i[$ (resp. $]-\i,0[$) is irregular.
Let us mention that the following inequalities have been proved 
in [F] :
$$\P(\int_{[0,+\i[} 1_{S_t\leq x} e^{-qt} L^s(dt)\leq 4 H_q(x)\quad\hbox{and}\quad  \P(\int_{[0,+\i[ }1_{-I_t\leq x} e^{-qt} L^i(dt)\leq 4 \cc H_q(x)$$

 \subsection{The main results}

\begin{theorem} \label {ED} For all $x\in ]0,+\i[$, $q\in [0,+\i[$, one has :

1) For all complex $\l$,  the functions  $A_q$ and  $\cc A_q$ satisfy the integral equations
$$A_q(x^-,\l) =H_q(0)+\int_{]0,x[}e^{-\l y}\cc A_q(y,\l){H_q(dy)\over \cc H_q(y)}$$
$$\cc A_q(x,\l) = \cc H_q(0)+\int_{]0,x]}e^{\l y} A_q(y^-,\l){\cc H_q(dy)\over H_q(y^-)}$$

2) For all complex  $\l$ with  $\Re(\l)>0$ (and  $\Re(\l)=0$ if  $q>0$ or $\lim X_t=-\i$) the functions  $C_q$ and $B_q$ satisfy the integral equations 

$$C_q(x^-,\l)=\int_{[x,+\i[}e^{-\l y}B_q(y,\l){ H_q(dy)\over \cc H_q(y)}$$
$$B_q(x,\l)= \psi_q(\l)+\int_{]x,+\i[} e^{\l y} C_q(y^-,\l){\cc H_q(dy)\over H_q(y^-)}$$

Moreover, one has 
$$A_q(x^-,\l)B_q(x,\l)+\cc A_q(x,\l) C_q(x^-,\l)=1$$

3) For  all complex $\l$ with $\Re(\l)<0$ (and $\Re(\l)=0$ if  $q>0$ or $\lim X_t=+\i$) the functions  $\cc C_q$ and $\cc B_q$ satisfy the integral equations 

$$\cc C_q(x,\l)=\int_{]x,+\i[}e^{\l y}\cc B_q(y^-,\l){\cc H_q(dy)\over H_q(y^-)}$$
$$\cc B_q(x^-,\l)= \cc \psi_q(\l)+\int_{[x,+\i[} e^{-\l y}\cc C_q(y,\l){H_q(dy)\over \cc H_q(y)}$$
Moreover, one has 
$$\cc A_q(x,\l)\cc B_q(x^-,\l)+A_q(x^-,\l)\cc C_q(x,\l)=1$$

\end{theorem}
For all $x\in ]0,+\i[$, $q\in [0,+\i[$ define : 

$$M_q(x,\l)=\left( \begin{array} {clcr} A_q(x^-,\l) & -C_q(x^-,\l) \\ \cc A_q(x,\l)  & B_q(x,\l)\end{array} \right)\hbox{ if } \quad \Re (\l)>0$$

$$M_q(x,\l)=\left( \begin{array} {clcr} \cc B_q(x^-,\l) & A_q(x^-,\l) \\ -\cc C_q(x,\l)  & \cc A_q(x,\l)\end{array} \right)\hbox{ if  } \quad \Re (\l)<0$$
 
Observe that for all  $\l=iu\in i{\bf R}$, the following limits exist : 
$$M_q^+(x,iu):=\lim_{\l \to iu, \Re(\l)>0}M(x,\l)=\left( \begin{array} {clcr} A_q(x^-,iu) & -C_q(x^-,iu) \\ \cc A_q(x,iu)  & B_q(x,iu)\end{array} \right)$$
$$ M_q^-(x,iu):= \lim_{\l \to iu, \Re(\l)<0}M(x,\l)=\left( \begin{array} {clcr} \cc B_q(x^-,iu) & A_q(x^-,iu) \\ -\cc C_q(x,iu)  & \cc A_q(x,iu )\end{array} \right) $$
The following result gives a Riemann-Hilbert characterization of the matrix $M_q$.

\begin{theorem}\label {RH} For all  $x\in ]0,+\i[$,  $q\in [0,+\i[$,  $\l\mapsto M_q(x,\l)$ is the unique function satisfying the following properties 

1) $\l\mapsto M_q(x,\l)$ is analytic on the two half-planes  $\{\Re(\l)>0\}$ and $\{\Re(\l)<0\}$.

2) $\l\mapsto M_q(x,\l)$  has a right limit ($M_q^+(x,iu)$)  and a left limit ($M_q^-(x,iu)$) at every point   $iu \in i{\bf R}$ and these two limits satisfy the equation : 

$$ M_q^+(x,iu) =M_q^-(x,iu) \left( \begin{array} {clcr} 0 & -1\\ 1& \phi(iu)+q \end{array} \right)$$

3) The matrix 
$$ \left( \begin{array} {clcr} {M_{11}\over |\l|+1 }& e^{\l x} M_{12}\\ e^{-\l x} M_{21}& {M_{22}\over |\l|+1 } \end{array} \right)$$
is bounded in $\{\Re(\l)>0\}\cup\{\Re(\l)<0\}$.

4) The following limits are valid for $\Re(\l) \to -\i$ :
$${ M_{11}(x,\l)\over \cc \psi_0(\l)} \to 1 \qquad \cc \psi_0(\l)M_{22}(x,\l)\to 1 \qquad e^{\l x} M_{12}\to 0$$
and $e^{-\l x} M_{21}(x,\l)\to 0$ if $]-\i, 0[$ is irregular.
\end{theorem}

\subsection{Connections with scattering theory} 

The integral equations of theorem \ref {ED} can be rewritten as the following distribution theoretic differential equation; for all  $\l \in \C\backslash i\R$,
$$M'_q(x,\l)= \left( \begin{array} {clcr} 0 & e^{-\l x} {H_q(dx)\over \cc H_q(x)} \\ e^{\l x} {\cc H_q(dx)\over  H_q(x^-)} &0 \end{array} \right)M_q(x,\l)\eqno{(Sc)}$$
This differential equation is a non standard form of the  classical equation of the scattering theory on the line,
with a   measure valued potential matrix $\left( \begin{array} {clcr} 0 & {H_q(dx)\over \cc H_q(x)} \\ {\cc H_q(dx)\over H_q(x^-)} &0 \end{array} \right)$.  Remark that these measures are unbounded in general.

Let us recall the basics of scattering theory, as expounded in \cite{BDZ}. 
One considers  a potential matrix $\left( \begin{array} {clcr} 0 & v \\ \cc v  &0 \end{array} \right)$ where $v$ et $\cc v$ are real functions of real variable satisfying some regularity assumptions, in particular they are integrable. To this matrix is associated the 
 differential equation :
$$Y'(x,\l)= \left( \begin{array} {clcr} 0 & e^{-\l x} v \\ e^{\l x}\cc v  &0 \end{array} \right)Y(x,\l)$$ where $\l$ is a complex parameter.

Denotes by  $J$ the matrix $J:=\left( \begin{array} {clcr} 1 & 0\\ 0 & -1 \end{array} \right)$. One can prove that  for any  imaginary complex number $\l=iu$ and  any matrix solution  $x\mapsto Y(x,iu)$, the matrix  $e^{{iu x\over 2} J} Y(x,iu) e^{{-iu x\over 2} J}$ 
converges as  $x\to+\i$ and $x\to-\i$ to the limits   $y^+(iu)$ and $y^-(iu)$. The scattering matrix associated to the potential is given by $S(iu)=y^+(iu) [y^-(iu)]^{-1}$.
Obviously this matrix doesn't depend on the particular solution $x\mapsto Y(x,\l)$.
On the other hand there exists a unique solution $x\mapsto {\mathcal Y}(x,\l)$ such that $e^{{\l x\over 2} J} {\mathcal Y}(x,\l) e^{{-\l x\over 2}J}$  converges to the identity matrix for $x\to-\i$.

Shabat \cite {Sh}, see also \cite{BDZ},  has proved  that for any  $x\in \R$, the function  $\l\mapsto {\mathcal Y}(x,\l)$ is analytic on the half planes $\{\Re(\l)>0\}$ and $\{\Re(\l)<0\}$, 
has right limits ${\mathcal Y}(x,(iu)^+)$ and left limits  ${\mathcal Y}(x,(iu)^-)$ at any point  $iu$ on the imaginary axis and the jumps matrix
 ${\mathcal Y}(x,(iu)^+)[{\mathcal Y}(x,(iu)^-)]^{-1}$ is equal to $S(iu)$. Moreover,   $e^{{\l x\over 2} J} {\mathcal Y}(x,\l) e^{{-\l x\over 2}J}$ converges  to the identity  matrix for  $|\l|\to +\i$. It is then clear that the function
$\l\mapsto{ \mathcal Y}(x,\l)$ is entirely determined by these properties. The determination of the matrix ${\mathcal Y}(x,\l)$ and consequently of the potential matrix from the scattering matrix ("the inverse problem ") is reduced to the solution of this Riemann-Hilbert problem.

In our problem, after extending  the potential by zero on $]-\i,0[$, one can see that the potential is not given by a function but rather by a measure, furthermore it is not integrable on $\R$ unless $\P$ is the distribution of a compound Poisson process and $q$ is positive. Consequently, none of the solutions of the associated differential equation is regular both for $x\to +\i$ and $x\to-\i$ and  one does not have any   scattering matrix in the classical sens (i.e. $S(iu)=y^+(iu) [y^-(iu)]^{-1}$). Moreover there is no  solution of the differential equation $x\mapsto Y(x,\l)$
such that the matrix $e^{{\l x\over 2} J} Y(x,\l) e^{{-\l x\over 2}J}$  converges to the identity matrix for $x\to-\i$.
We have chosen the solution of the differential equation  $(Sc)$ that  is  the most convenient from a probabilistic point of view.
Theorem \ref {RH}  tells us that the determination of this solution (and consequently of all the others and of the potential matrix) is reduced to the solution of a Riemann-Hilbert problem, as in classical scattering theory. In our setting 
the role of the scattering matrix is played by the matrix $\left( \begin{array} {clcr} 0 & -1\\1 & \phi(iu ) +q\end{array} \right)$.
\section{Connection with the  exit time from an interval and related variables}

Let us  now explain how  distributions of certain random variables related to the exit time of a L\'evy process from an interval are related to the matrix $M_q(x,\l)$.

Let us denote 
$\hbox{supp}( L^s)$ and  $\hbox{supp}( L^i)$ the supports of the random measures $L^s(dt)$ and $L^i(dt)$, and for all positive time  $t$, 
$$g^i_t=\sup\left([0,t[\cap \hbox{supp} (L^i)\right)\qquad g^s_t=\sup\left([0,t[\cap \hbox{supp}( L^s)\right)$$

\begin{proposition}\label{PT} For all $\l_1,\l_2 \in {\bf C}$,  $q_1, q_2\in [0,+\i[$, $x\in ]0,+\i[$, one has 

$$\cc A_{q_1}(x,\l_1)A_{q_2}(x,\l_2)=\P( \int _{[0,+\i [}e^{-\l_1 I_t-q_1g^i_t} 1_{S_t-I_t\leq x} e^{-\l_2 (X_t-I_t)-q_2(t-g^i_t)} dt)$$
$$=\P( \int _{[0,+\i[}e^{-\l_2S_t-q_2 g^s_t}1_{S_t-I_t\leq x} e^{-\l_1 (X_t-S_t)-q_1(t-g^s_t )} dt)$$
 \end{proposition} 
 This proposition indicates in particular, for $q_1=q_2=q$ that the pair of functions $(A_q,\cc A_q)$ determines the resolvent of the trivariate Markov process $(S,X,I)$. We shall see that  this result is an immediate consequence of a path decomposition stated in   proposition \ref {DT} and we will prove it later.
 
The  next propositions give expressions of  the distribution of some random variables related to fluctuations of the L\'evy process  in terms of our functions $A,B,C,\cc A,\cc B,\cc C$. They are immediate consequences of the definition of these functions and of the compensation formula stated in proposition \ref {CF}. We leave the proofs to the reader.
 
  \begin{proposition}For all $\Re(\mu_1)\geq 0$, $\Re(\mu_2)\leq 0$, $q_1,q_2\in [0,+\i[$, $x\in ]0,+\i[$, one has 
  
 $$\P\Bigr(e^{-\mu_1S_{T^s_x}-q_1 g^s_{T^s_x}}e^{-\mu_2( X_{T^s_x}-S_{T^s_x}) -q_2 (T^s_x- g^s_{T^s_x})}\Bigr)={1\over B_{q_1}(x,\mu_1)} \cc C_{q_2}(x,\mu_2)$$
 
  $$\P\Bigl(e^{-\mu_2 I_{T_i^x}-q_1 g^i_{T_i^x} }e^{-\mu_1( X_{T_i^x}-I_{T_i^x})-q_2 (T_i^x- g^i_{T_i ^x})}\Bigr)={1\over \cc B_{q_1}(x,\mu_2)}C_{q_2}(x,\mu_1)$$
 \end{proposition}

Let us denote 
$$U^x:=\inf\{t;  S_t-I_t\geq x \hbox { and  }S_t=X_t \}\wedge \inf\{t;  S_t-I_t>x \hbox { and }X_t=I_t \}$$

\begin{proposition} \label {A}
For all  $\l\in\C$, $\Re(\mu_1)\geq 0$,  $\Re(\mu_2)\leq 0$, $q_1,q_2\in [0,+\i[$, $x\in ]0,+\i[$, 
$$\P\Bigl( e^{-\l I_{U^x}-q_1 g^i_{U^x}}.e^{-\mu_1 (X_{U^x}- I_{U^x})-q_2 (U^x-g^i_{U^x})}; S_{U^x}=X_{U^x}\Bigr )= \cc A_{q_1}(x,\l) C_{q_2}(x^-,\mu_1)$$
$$\P\Bigl ( e^{-\l S_{U_x}-q_1 g^s_{U^x} } e^{-\mu_2 (X_{U^x}- S_{U^x})-q_2 (U^x-g^s_{U^x})}; I_{U^x}=X_{U^x}\Bigr)=A_{q_1}(x^-,\l) \cc C_{q_2}(x,\mu_2)$$
\end{proposition}

Let $A_q(x,dy)$ be the measure with Laplace transform 
$A_q(x,\l)$ ($ A_q(x,dy):=\P(\int_0^{+\i}1_{S_t-I_t\leq x}1_{S_t\in dy}L^s(dt))$). This measure is absolutely continuous with respect to  the measure $G_q(dy):=\P( \int _{[0,+\i[}1_{S_t\in dy}e^{-qt} L^s(dt))$, with a density $\alpha_q (x,y)$, this entitles us to define

$$A_q(b+y,dy):=\alpha_q (b+y,y)G_q(dy)=\P( \int _{[0,+\i[}1_{-I^-_t\leq b} 1_{S_t\in dy} e^{-qt} L^s(dt))$$
Similarly, denote : 
$\cc A_q(a-y,dy)$ the measure :
 $$\cc A_q(a-y,dy):=\P( \int _{[0,+\i[}1_{S_t\leq a} 1_{I_t\in dy} e^{-qt}L^i(dt))$$

For all  positive reals $a,b$, denote $T_b^{a}:=\inf\{t; X_t\not\in [-b,a]\}$.

\begin{proposition}  For all  $\l\in\C$, $\Re(\mu_1)\geq 0$,  $\Re(\mu_2)\leq 0$, $q_1,q_2\in [0,+\i[$

$$\P\Bigl(e^{-\l I_{T_b^a}-q_1g^i_{T_b^a}}. e^{-\mu_1 (X_{T_b^a}-I_{T_b^a})-q_2({T_b^a} -g^i_{T_b^a})};X_{T_b^a}=S_{T_b^a})= \int_{[-b,0]} e^{-\l y} C_{q_2}(a-y,\mu_1)\cc A_{q_1}(a-y,dy)$$

$$\P\Bigl(e^{-\l S_{T_b^a}-q_1g^s_{T_b^a}}. e^{-\mu_2 (X_{T_b^a}-S_{T_b^a})-q_2({T_b^a} -g^s_{T_b^a})};X_{T_b^a}=I_{T_b^a})= \int_{[0,a]} e^{-\l y} \cc C_{q_2}(b+y,\mu_2)A_{q_1}(b+y,dy)$$

\end{proposition}

\section{Independence of past and post minimum process "conditionnally on the amplitude" and  
proof of proposition \ref{PT}}

Let us first introduce some more notations : For all $q\in [0,+\i[$,  let
 $\Q_q$, $\Q^{\uparrow}_q$ and $\Q^{\downarrow}_q$  be the measures on $(\Omega, {\mathcal F},\P)$ defined as follows :
 $$\Q_q(dw):=\P(\int_0^{+\i}1_{ (X_s)_{0\leq s<t}\in dw} e^{-qt}dt)$$
$$\Q^{\uparrow}_q(dw):= \P[\int _{[0,+\i[} 1_{(S_t-X_{(t-s)^-}; s\geq 0)\in dw}e^{-qt}L^s(dt)]$$ 
 $$\Q^{\downarrow}_q(dw):= \P[\int _{[0,+\i[} 1_{(I_t-X_{(t-s)^-}; s\geq 0)\in dw}e^{-qt}L^i(dt)]$$
 
All these measures are supported by the set of paths with finite life time. 
Note  that for a positive $q$, $\P_q:=q\Q_q$ is 
 the distribution of the L\'evy process $X$ under $\P$ killed at an independent exponential time with parameter $q$ and  that the measure $\Q^{\uparrow}_0$ is finite if  $\lim_{t\to+\i}X_t=-\i$ ($\P$-a.s.), and $\Q^{\downarrow}_0$ is finite if  $\lim_{t\to+\i}X_t=+\i$ ($\P$-a.s.). The measures $\Q^{\uparrow}_0$ and $\Q^{\downarrow}_0$ are infinite  in the other cases and $\Q_0$ is infinite in all cases.

Denote $\z$ the life time of the canonical process $X$, and  $F$, $M$, $m$ the final values of $X$, $S$, $I$ :
$$\z:=\sup\{ t, X_t\in \R\}\qquad F:=X_{\z-}\qquad M:=S_{\z-}\qquad m:=I_{\z-}$$
Denote $\sigma$ and $\rho$ respectively the last time $X$  takes its maximal value and the first time $X$   takes its minimal  value : 
$$\sigma:=\sup\{ t\leq \z; X_t\vee X_{t-}=M\}\qquad \rho:= \inf\{ t\geq 0; X_t\wedge X_{t-}=m\}$$
We make the convention that $\z,\sigma,\rho, F, M, m$ are zero when the path is constantly equal the cemetery point $\delta$. 

Note that under $\Q_q$  for every $q\in [0,+\i[$,  the time $\sigma$ (respectively $\rho$) is respectively the unique time at which the process $X_t\vee X_{t-}$ (respectively $X_t\wedge X_{t-}$), takes its maximal value, (respectively its minimal value) unless when $\P$ is the distribution of a compound Poisson process. Note also that when $q$ is positive, the measures 
${\Q_q^{\uparrow}(dw)\over \Q_q^{\uparrow}(\Omega)}$ and ${\Q_q^{\downarrow}(dw)\over \Q_q^{\downarrow}(\Omega)}$ are respectively the law of
the process $ (M-X_{(\sigma-s)^-}, s\geq 0)$ and  the law of $ (m-X_{(\rho-s)^-}; s\geq 0)$ under $\P_q=q\Q_q$. This last property is a standard consequence of the compensation formula of proposition \ref{CF}.

The next identities follow  directly from the definitions.

$$A_q(x,\l)=\Q_q^{\uparrow}(e^{-\l F};M\leq x )=\Q_0^{\uparrow}(e^{-\l F-q\z}; M\leq x)$$
$$H_q(x)=\Q_q^{\uparrow}(M\leq x)=\Q_0^{\uparrow}(e^{-q \z};M\leq x)$$
$${1\over \psi_q(\l)}=\Q_q^{\uparrow}(e^{-\l F})=\Q_0^{\uparrow}(e^{-\l F-q\z}) $$
$$\cc A_q(x,\l)=\Q_q^{\downarrow}(e^{-\l F};-m\leq x)=\Q_0^{\downarrow}(e^{-\l F-q\z};-m\leq x)$$
$${1\over \cc \psi_q(\l)}=\Q_q^{\downarrow}(e^{-\l F})=\Q_0^{\downarrow}(e^{-\l F-q\z}) $$
$$ \cc H_q(x)=\Q_q^{\downarrow}(-m\leq x)=\Q_0^{\downarrow}(e^{-q\z};-m\leq x)$$

\begin{proposition} \label {DT} For every $q\in [0,+\i[$, $x\in ]0,+\i]$, 
$$\Q_q( (m-X_{(\rho-t)^-}; t\geq 0) \in dw_1; M-m\leq x, (X_{t+\rho}-m; t\geq 0)\in dw_2)$$
$$=\Q_q( (M-X_{(\sigma-t)^-}; t\geq 0) \in dw_2; M-m\leq x, (X_{t+\sigma}-M; t\geq 0)\in dw_1)$$
$$=\Q_{q}^{\downarrow}(dw_1;-m\leq x )\Q_{q}^{\uparrow}(dw_2; M\leq x)$$
\end{proposition}
{\bf Proof } Let us first prove the assertions for $q>0$ and $x=+\i$.
One has the following identities : 
$$\P_q( (m-X_{(\rho-t)^-}; t\geq 0) \in dw_1; (X_{t+\rho}-m; t\geq 0)\in dw_2)$$
$$=\P_q( (m-X_{(\rho-t)^-}; t\geq 0) \in dw_1)
\P_q((X_{t+\rho}-m; t\geq 0)\in dw_2)$$
$$=\P_q( (m-X_{(\rho-t)^-}; t\geq 0) \in dw_1)
\P_q((M-X_{(\sigma-t)^-}, t\geq 0)\in dw_2)$$
$$={\Q_q^{\downarrow}(dw_1)\over \Q_q^{\downarrow}(\Omega)}{\Q^{\uparrow}_q(dw_2)\over \Q^{\uparrow}_q(\Omega)}$$ 

The first equality follows from the well known independence of the past and post minimum processes (see [B] lemma 6 chapter 6 for exemple), the second from the fact,  again
 well known (see [B] lemma 2 chapter 2), that the process $(F-X_{(\z-t)^-}; t\geq 0)$ has the same law as $X$ under $\P_q$.
  The last one follows from the compensation formula, as we have already mentioned.

On the other hand, one gets 

$$\Q^{\downarrow}_q(\Omega) \Q^{\uparrow}_q(\Omega)={1\over \cc \psi_q(0))}{1\over \psi_q(0)}={1\over q}$$
The first identity follows from the definitions
 of $\Q^{\downarrow}_q$ and $\Q^{\uparrow}_q$, the second one from   the Wiener-Hopf equation  of proposition \ref{WH} for $iu=0$.

Remember that $\P_q=q\Q_q$ and simplify the previous identities by ${1\over q}$ to get 
$$\Q_q( (m-X_{(\rho-t)^-}; t\geq 0) \in dw_1; (X_{t+\rho}-m; t\geq 0)\in dw_2)=\Q_q^{\downarrow}(dw_1)\Q^{\uparrow}_q(dw_2)$$

Letting $q$ goes to 0, we get the same identity for $q=0$. 

The events $\{M-m\leq x\}$ can be written as the intersection :
$$\{M-m\leq x\}=\{-\min [m-X_{(\rho-t)^-}; t\geq 0]\leq x\}\cap \{\max[X_{t+\rho}-m; t\geq 0]\leq x\}$$
this yields

$$\Q_q( (m-X_{(\rho-t)^-}; t\geq 0) \in dw_1; M-m\leq x; (X_{t+\rho}-m; t\geq 0)\in dw_2)$$
$$=\Q_q^{\downarrow}(dw_1;-m\leq x)\Q^{\uparrow}_q(dw_2;M\leq x)$$

Using the identity in law of the process $(F-X_{(\z-t)^-}; t\geq 0)$ and $X$ one gets the other identity of the proposition \qed

{\bf In the sequel, unless explicitely mentioned, all the properties hold for every  non negative $q$ and we shall omit to mention it}

{\bf Proof  of proposition \ref {PT} }

Denote $X^{\uparrow}$ and $X^{\downarrow}$  the processes $(X_{t+\rho}-m; t\geq 0)$ and $ (m-X_{(\rho-t)^-}; t\geq 0)$ respectively. One gets 
$$\P( \int _{[0,+\i [}e^{-\l_1 I_t-q_1g^i_t} 1_{S_t-I_t\leq x} e^{-\l_2 (X_t-I_t)-q_2(t-g^i_t)} dt)$$
$$=\Q_0( e^{-\l_1 m-q_1\rho}1_{M-m\leq x} e^{-\l_2(F-m)-q_2(\z-\rho)})$$
$$=\Q_0(e^{-\l_1F(X^{\downarrow})-q_1\z(X^{\downarrow})} 1_{-m(X^{\downarrow})\leq x} 1_{M(X^{\uparrow})\leq x}e^{-\l_2F(X^{\uparrow})-q_2\z(X^{\uparrow})})$$
$$=\Q_0^{\downarrow}( e^{-\l_1 F-q_1\z};-m\leq x)\Q_0^{\uparrow}( e^{-\l_2 F-q_2\z}; M\leq x)=\cc A_{q_1}(x,\l_1)A_{q_2}(x,\l_2)$$
The first, second and fourth equalities follow from the definitions. The third one follows from proposition \ref{DT}.

Similarly, the second assertion of proposition \ref{PT} follows from the second assertion of proposition \ref{DT}.
\qed

\section{\bf Proof of property  1) of theorem \ref {ED}}

Denote for any $x\in ]0,+\i[$, 
$$\Q_q^{\downarrow x}:=\Q_q^{\downarrow}(dw\vert -m\leq x)\qquad \Q_q^{\uparrow x}:=\Q_q^{\uparrow}(dw\vert M<x)$$

Note that when $\Q_q^{\downarrow}(m=0)>0$, the mesure  $\Q_q^{\downarrow}(dw\vert -m=0)$ is the Dirac mass on the path constantly equal to $\delta$. We make the convention that $\Q_q^{\downarrow 0}$ and  $\Q_q^{\uparrow 0}$ are this Dirac mass.

\begin{lemma}\label {CD}

$$\Q_q^{\uparrow}((X_{\sigma+t}-M; t\geq 0)\in dw\vert M-X_{(\sigma-t)^-}; t\geq 0)=\Q_q^{\downarrow M}(dw)$$
 
$$\Q_q^{\downarrow}((X_{\rho+t}-m; t\geq 0)\in dw\vert m-X_{(\rho-t)^-};  t\geq 0)=\Q_q^{\uparrow -m}(dw)$$

\end{lemma}

{\bf Proof } Let us discuss few facts about the event  $\{\rho\leq \sigma\}$. First of all, one has 
$$\{\rho\leq \sigma\}= \{I_{\sigma}\leq  \inf (X_s; s\geq \sigma)\}$$ and 
$$\{\rho\leq \sigma\}=\{S_{\rho}^-\leq \sup (X_s;s\geq \rho)\vee X_{\rho^-}\}$$

If $]-\i,0[$ is regular then $\Q_q$-a.s. the canonical process $X$ has no negative jump at the times $\rho$ and $\sigma$.
Consequently,
$I_{\sigma}=I_{\sigma^-}$ and $\sup \{X_s;s\geq \rho\}\vee X_{\rho^-}=\sup \{X_s;s\geq \rho\}$ and one gets :
$$\{\rho\leq \sigma\}= \{I_{\sigma}^-\leq  \inf (X_s; s\geq \sigma)\}=\{S_{\rho}^-\leq \sup (X_s;s\geq \rho)\}$$

If  $]-\i,0[$ is irregular then, $\Q_q$-a.s., either $\rho$ is zero either it is a time when $X$ has a negative jump and $\sigma$ is  a time when $X$ has a negative jump. On the other hand, when $\rho=\sigma$, we have $I_{\sigma^-}> \inf (X_s; s\geq \sigma)$ and  $S_{\rho}^-> \sup (X_s;s\geq \rho)$. Thus

$$\{\rho<\sigma\}= \{I_{\sigma}^-\leq  \inf (X_s; s\geq \sigma)\}=\{S_{\rho}^-\leq \sup (X_s;s\geq \rho)\}$$

For the needs of  next  proof, let us denote $A$ the event 
$$A:=\{I_{\sigma}^-\leq  \inf (X_s; s\geq \sigma)\}=\{S_{\rho}^-\leq \sup (X_s;s\geq \rho)\}$$
By the preceeding discussion, the event $A$ is either  $\{\rho\leq \sigma\}$ or $\{\rho<\sigma\}$.

Let  $X^{\downarrow}$ and $X^{\uparrow}$  denote respectively 
the processes $(X_{\sigma+t}-M,t\geq 0)$ and  $(X_{\rho+t}-m; t\geq 0)$, we add a $^{\downarrow}$ or $^{\uparrow} $ to the corresponding objects.

One has $$A=\{ I_{\sigma^-}-M\leq m(X^{\downarrow})\}=\{S_{\rho}^--m\leq M(X^{\uparrow})\}$$
Denote   $\overline {\mathcal F}_{\sigma}$  the $\sigma$-field generated  by the pre-maximum process  $(M-X_{(\sigma-t)^-}, t\geq 0)$.
Notice that the trace on the event $A$ of the $\sigma$-field $\overline {\mathcal F}_{\sigma}$ contains the random variable $M-m1_A$ and consequently, this trace is $\sigma$-finite for $\Q_q$.
The independence under $\Q_q$ of the process $X^{\downarrow}$  and the $\sigma$-field
$\overline{\mathcal F_{\sigma}}$ given in proposition \ref {DT} gives us :
 $$\Q_q ( X^{\downarrow}\in dw \vert \overline{\mathcal F}_{\sigma}\vee \sigma(A))=\Q_q ( X^{\downarrow}\in dw \vert \overline{\mathcal F}_{\sigma}\vee \sigma(\{ I_{\sigma^-}-M\leq m(X^{\downarrow})\})= \Q_q^{\downarrow M- I_{\sigma^-}}(dw)\quad \hbox 
 {on }\quad  A$$
 Moreover, the trace on $A$ of the $\sigma$-field $\overline {\mathcal F}_{\sigma}$ contains the $\sigma$-field $\overline{\mathcal F}_{\sigma}^{\uparrow}\vee \sigma (S_{\rho}^--m) $, and 
 $$(X^{\uparrow})^{\downarrow}=X^{\downarrow} \quad \hbox{and}\quad M-I_{\sigma^-}=M(X^{\uparrow})\qquad \hbox {on } A$$
 
The previous identity gives then the following :

$$\Q_q( (X^{\uparrow})^{\downarrow}\in dw \vert \overline {\mathcal F}_{\sigma}^{\uparrow}\vee \sigma( S_{\rho}^--m))=
\Q_q^{\downarrow, M(X^{\uparrow})}(dw)\quad 
 \hbox{ on  } A=\{M(X^{\uparrow})\geq S_{\rho}^--m\} $$
 
Using the independence of the process   $X^{\uparrow}$ and the random variable  $S_{\rho}^--m$ given by proposition \ref {DT}, one then deduces easily the first identity  of the lemma.

One gets the second identity similarly with using $^cA$ instead of A.\qed

{\bf Proof of property  1) of  theorem \ref {ED}}
For all  $x\in ]0,+\i[$ and  $\l\in \C$, one has

$$
\begin{array}{lrc}
A_q(x^-,\l)&=&\Q^{\uparrow}_q(M<x,e^{-\l F})=H_q(0)+\int_{]0,x[}\Q^{\uparrow}_q(e^{-\l F}\vert M=y)H_q(dy)\\
&=&H_q(0)+\int_{]0,x[}e^{-\l y}\Q^{\uparrow}_q(e^{-\l (F-M)}\vert M=y)H_q(dy)\\
&=&H_q(0)+\int_{]0,x[}e^{-\l y}\Q_q^{\uparrow}(e^{-\l F[X_{\sigma +t}-M;t\geq 0]}\vert M=y)H_q(dy)
\end{array}
$$
Using Lemma \ref {CD}, this last quantity is equal to
$$
\begin{array}{lrc}
&=&H_q(0)+\int_{]0,x[}e^{-\l y}\Q_q^{\downarrow}(e^{-\l F}\vert -m\leq y)H_q(dy)=H_q(0)+\int_{]0,x[}e^{-\l y}\Q_q^{\downarrow}(e^{-\l F};-m\leq y){H_q(dy)\over \cc H_q(y)}\\
&=&H_q(0)+\int_{]0,x[}e^{-\l y}\cc A_q(y,\l){H_q(dy)\over \cc H_q(y)}
\end{array}
$$

We get the identity $$\cc A_q(x,\l)=\cc H_q(0)+\int_{]0,x]}e^{\l y} A_q(y^-,\l){\cc H_q(dy)\over H_q(y^-)}$$ in a similar way.

\qed

\section {Proof of properties 2 and 3 of theorem \ref {ED}} 
\subsection{Two  Markov chains}
We now define the successive minima and maxima.
We let first 
$$
\begin{array}{rcl}
M_1:=M&\  \hbox{and}&\ T_1:=\sigma\\
m_2:=\inf\{X_t; t\geq T_1\}&\ \hbox{and}\ & T_2:=\inf\{t\geq T_1; X_t\wedge X_{t-}=m_2\}
\end{array}
$$ then we define inductively
$$
\begin{array}{rcl}
M_{2n+1}:=\sup\{X_t; t\geq T_{2n}\}&\ \hbox{and}&\ T_{2n+1}:=\sup\{t\geq T_{2n}; X_t\vee X_{t-}=M_{2n+1}\}\\
m_{2n+2}:=\inf\{X_t; t\geq T_{2n+1}\}&\ \hbox{and}&\ T_{2n+2}:=\inf\{t\geq T_{2n+1}; X_t\wedge X_{t-}=m_{2n+2}\}
\end{array}
$$ 
and $$Z_{2n+1}:=M_{2n+1}-m_{2n} \quad Z_{2n+2}:=m_{2n+2}-M_{2n+1}$$
If $T_n=\z$ put  $Z_{n+1}:=0$ and  $T_{n+1}:=\z$.
Notice that $F=\sum_1^{+\i} Z_n$.

Below is picture of our sequence.

$$\includegraphics[scale=0.8]{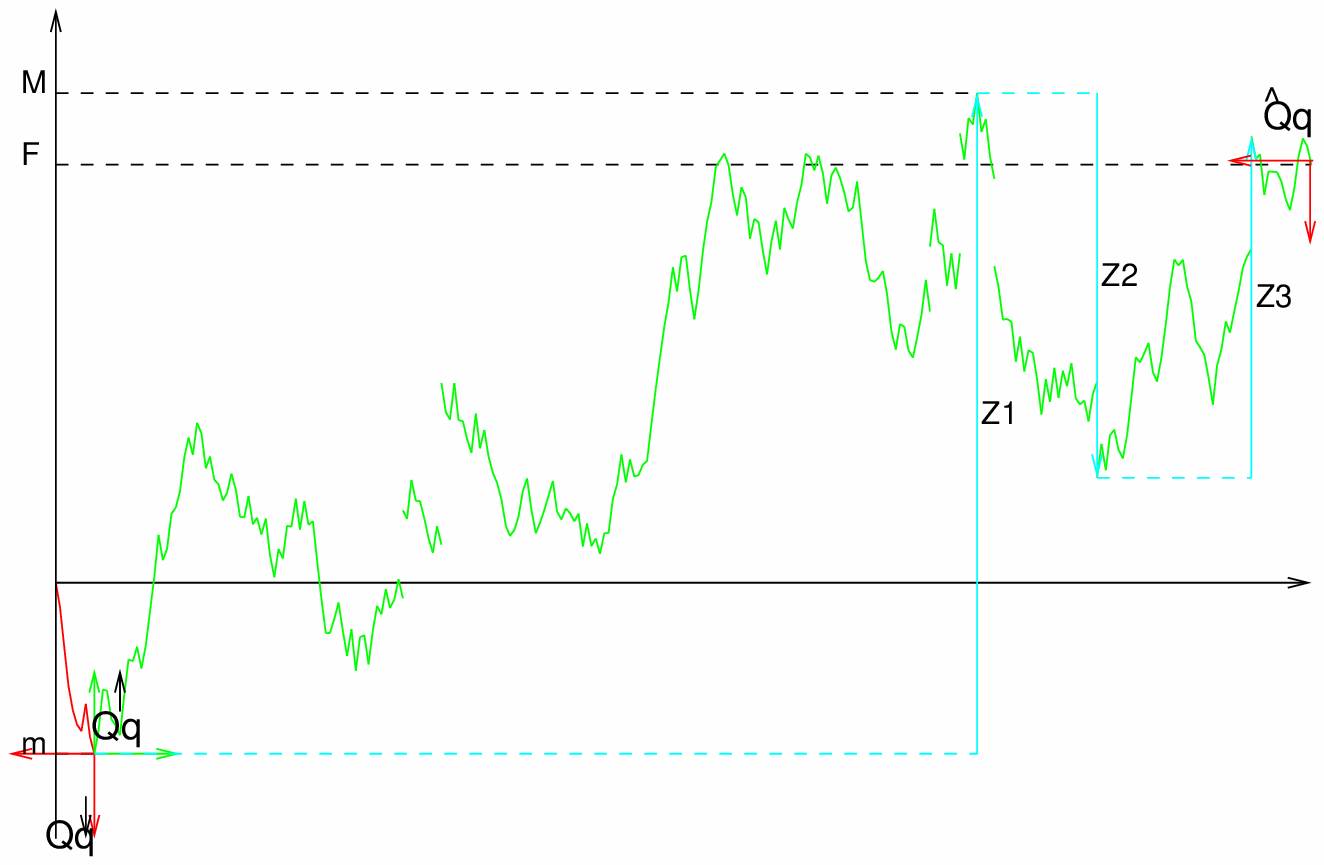}$$
Similarly, let
$$
\begin{array}{rcl}
m_1:=m&\  \hbox{and}&\ S_1:=\rho\\
M_2:=\sup\{X_t; t\geq T_1\}&\ \hbox{and}\ & S_2:=\sup\{t\geq T_1; X_t\vee X_{t-}=I_2\}
\\
m_{2n+1}:=\inf\{X_t; t\geq S_{2n}\}&\ \hbox{and}&\ S_{2n+1}:=\inf\{t\geq S_{2n}; X_t\wedge X_{t-}=m_{2n+1}\}\\ 
M_{2n+2}:=\sup\{X_t; t\geq S_{2n+1}\}&\ \hbox{and}&\ S_{2n+2}:=\sup\{t\geq S_{2n+1}; X_t\vee X_{t-}=M_{2n+2}\}

\end{array}
$$ 
and 
 $$Y_{2n+1}:=m_{2n+1}-M_{2n}\quad Y_{2n+2}:=M_{2n+2}-m_{2n+1}$$
 If $S_n=\z$ put  $Y_{n+1}:=0$ and  $S_{n+1}:=\z$. Notice that $F=\sum_1^{+\i} Y_n$.

\begin{lemma} \label {MC} Under  $\Q_q^{\uparrow}$ 
the sequence  $(Z_1,Z_2,\dots)$ is a Markov chain
with  transition kernel given by 
$$\hbox{for  $x\not=0$}\quad P_q(x,dy)={H_{q}(dy)\over H_{q}((-x)^-)}1_{y\in [0,-x[} + {\cc H_{q}(-dy)\over \cc H_{q}(x)} 1_{y\in [-x,0]}$$
$$P_q(0,\{0\})=1$$
and with  initial  law : $$\Q_q^{\uparrow}(Z_1\in dy)=H_q(dy)$$
Under  $\Q_q^{\downarrow}$ the sequence $(Y_1,Y_2,\dots )$ is a  Markov chain with the same transition kernel
and with initial   law $$\Q_q^{\downarrow}(Y_1\in dy)=\cc H_q(dy)$$
\end{lemma}
{\bf Proof } One has : 
$$\Q_q^{\uparrow}(Z_1\in dy)= H_q(dy)$$
$$
\begin{array}{c}
\Q_q^{\uparrow}( Z_2\in dy\vert Z_1=x)=\Q_q^{\uparrow}(\min[X_{\sigma+t}-M; t\geq 0]\in dy \vert M=x)\\
=\Q_q^{\downarrow}(m\in dy\vert -m\leq x)={\cc H_q(-dy)\over \cc H_q(x)}1_{y\in [-x,0]}
\end{array}
$$
The third identity comes from Lemma \ref{CD}, the others follow from definitions of the random variables $Z_1$ and $Z_2$.

The same arguments yields : 
$$\Q_q^{\downarrow}( Y_2\in dy\vert Y_1=x)={H_q(dy)\over H_q(x^-)}1_{y\in [0,x[}$$
A straightforward induction  gives us the rest of the lemma.\qed

Let 
$$U_q(dy):=\sum_1^{+\i}\Q_q^{\uparrow}(Z_n\in dy)\qquad V_q(dy):=\sum_1^{+\i}\Q_q^{\downarrow}(Y_n\in dy)$$
 $$\tau_x:=\sup\{n; Z_n\not\in [-x,x[\}\qquad( \tau_x:=0 \hbox{ if }Z_1<x)$$ 
$$\nu_x:=\sup\{n; Y_n\not\in[-x,x[\}\qquad (\nu_x:=0\hbox{ if }Y_1\geq -x)$$
\begin{lemma} \label {MR}
Under $\Q_q^{\uparrow}(dw; \tau_x>0)$ for every $x\in ]0,+\i[$,  the sequence $(Z_{\tau_x-n})_{0\leq n<\tau_x}$ is a sub-Markov chain with values in 
 $\R^*$ and with  initial law :
 $$\Q_q^{\uparrow}(Z_{\tau_x}\in dy;\tau_x>0)=[1_{y\in ]-\i,-x[}{H_q(x^-)\over H_q((-y)^-)}+1_{y\in [x,+\i[}{\cc H_q(x)\over \cc H_q(y)}]U_q(dy)$$
Its transition kernel  $R_q$ does not depend on $x$ and satisfies the equation  :
$$U_q(dz) P_q(z,dy)=R_q(y,dz)U_q(dy)\qquad \hbox{ on }\quad \R^*\times\R^*$$

Under $\Q^{\downarrow}_q(dw,\nu_x>0)$, the sequence  $(Y_{\nu_x-n})_{0\leq n<\nu_x}$ is a sub-Markov chain 
with values in  $\R^*$ and with  initial  law :
$$\Q_q^{\downarrow}(Y_{\nu_x}\in dy,\nu_x>0)=[1_{y\in ]-\i,-x[}{H_q(x^-)\over H_q((-y)^-)}+1_{y\in [x,+\i[}{\cc H_q(x)\over \cc H_q(y)}]V_q(dy)$$
Its transition kernel  $S_q$ does not depend on $x$ and satisfies the equation  :
$$V_{q}(dz) P_{q}(z,dy)=S_{q}(y,dz)V_{q}(dy)\qquad \hbox{ on }\quad \R^*\times\R^*$$
\end{lemma}

{\bf Proof} It is a standard fact from the theory of time reversal of Markov chains (or processes) that the time 
$\tau_x$ is a so called "return time" of the Markov chain  $(Z_n)$ and consequently  that 
$(Z_{\tau_x-n})_{0\leq n<\tau_x}$ is a sub-Markov chain with a  transition kernel related to the  one of $(Z_n)$ by the so called "duality identity"
$$U_q(dz) P_q(z,dy)=R_q(y,dz)U_q(dy)\qquad \hbox{ on }\quad \R^*\times\R^*$$
In particular, this transition kernel does not  depend on the particular return time $\tau_x$. 
Since we lack  an adequate reference we check this property in our particular case. 

For every $y\in \R^*$, let 
 $z\mapsto p_q(y,z)$ be a density of the measure $P_q(y,dz)$ relatively to $U_q(dz)$, one has for every integer $m$ (measures involved  here are on $\R^*$) $$\Q^{\uparrow}_q( Z_{\tau_x-m}\in dy_m , Z_{\tau_x-(m-1)}\in dy_{m-1},\dots, Z_{\tau_x}\in dy_0)$$
 $$=\sum_{n=m+1}^{+\i}\Q^{\uparrow}_q( Z_{n-m}\in dy_m , Z_{n-(m-1)}\in dy_{m-1},Z_{n-(m-2)}\in dy_{m-2}, \dots, Z_{n}\in dy_0,\tau_x=n)$$
 $$=\sum_{n=m+1}^{+\i}\Q^{\uparrow}_q( Z_{n-m}\in dy_m , Z_{n-(m-1)}\in dy_{m-1},Z_{n-(m-2)}\in dy_{m-2},\dots, Z_{n}\in dy_0,Z_{n+1}\in [-x,x[)1_{y_0\not\in [-x,x[}$$
 $$=\sum_{n=m+1}^{+\i}\Q^{\uparrow}_q( Z_{n-m}\in dy_m) p_q(y_m,y_{m-1})U_q(dy_{m-1})
 P_q(y_{m-1},dy_{m-2})\dots P_q(y_1, dy_0)P_q(y_0,[-x,x[)1_{y_0\not\in [-x,x[}$$
 $$= U_q(dy_m)p_q(y_m,y_{m-1})U_q(dy_{m-1})
 P_q(y_{m-1},dy_{m-2})\dots P_q(y_1, dy_0)P_q(y_0,[-x,x[)1_{y_0\not\in [-x,x[} $$
 
 [The third equality follows from the definition of our particular  time $\tau_x$ and from the fact that the sequence $\vert Z_n\vert$ is nonincreasing].
 
 One can deduce from the last equalities that :
$$\Q^{\uparrow}_q( Z_{\tau_x-m}\in dy_m| Z_{\tau_x-(m-1)},\dots,Z_{\tau_x})=U_q(dy_m)p_q(y_m, Z_{\tau_x-(m-1)})  $$
Thus the sequence $(Z_{\tau_x-n})_{0\leq n<\tau_x}$ is a sub-Markov chain and its transition kernel is $R_q(z,dy)=p_q(y,z)U_q(dy)$. Multiplying the last identity by the measure $U_q(dz)$, one recognizes the duality equation.

Let us  compute the  law of $(Z_{\tau_x-n})_{0\leq n<\tau_x}$. Put $m=0$ in the preceeding equalities and get 
$$\Q_q^{\uparrow}(Z_{\tau_x}\in dy_0)=1_{y_0\notin [-x,x[}\sum_{n=1}^{+\i}\Q_q^{\uparrow} (Z_n\in dy_0; Z_{n+1}\in [-x,x[)$$
$$=1_{y_0\not\in [-x,x[} \sum_{n=1}^{+\i}\Q_q^{\uparrow}( Z_n\in dy_0) P_q(y_0, [-x,x[)$$
 $$= \Bigl(1_{y_0\in [x,+\i[}{\cc H_q(x)\over \cc H_q(y_0)}+1_{y_0\in ]-\i,-x[}{H_q(x^-)\over H_q((-y_0)^-)}\Bigr)U_q(dy_0)$$

The last identity comes from  the value of $P_q$ given by lemma \ref {MC}. One gets the results for the sequence  $(Y_{\nu_x-n})_{0\leq n<\nu_x}$ similarely.\qed
 
Let 
$$c_q(x,\l):={1\over \cc H_q(x)}\Q^{\uparrow}_q(e^{-\l\sum_{n=1}^{\tau_x} Z_n};\tau_x \hbox{is odd})={1\over \cc H_q(x)}\Q^{\uparrow}_q( e^{-\l M_{\tau_x}};\tau_x \hbox{is odd}) $$
$$d_q(x,\l):={1\over H_q(x^-)}\Q^{\uparrow}_q(e^{-\l\sum_{n=1}^{\tau_x} Z_n};\tau_x \hbox{is even})={1\over H_q(x^-)}\Q_q^{\uparrow}( e^{-\l m_{\tau_x}};\tau_x \hbox{is even}) $$

$$\cc c_q(x,\l):={1\over H_q(x^-)}\Q^{\downarrow}_q(e^{-\l\sum_{n=1}^{\nu_x} Y_n};\nu_x \hbox{is odd})={1\over H_q(x^-)}\Q^{\downarrow}_q( e^{-\l m_{\nu_x} };\nu_x \hbox{is odd}) $$
$$\cc d_q(x,\l):={1\over \cc H_q(x)}Q^{\downarrow}_q(e^{-\l\sum_{n=1}^{\nu_x} Y_n};\nu_x \hbox{is even})={1\over \cc H_q(x)}\Q^{\downarrow}_q( e^{-\l M_{\nu_x}};\nu_x \hbox{is even}) $$

\begin{lemma} \label {bc}For $\Re(\l) >0$ (and  $\Re(\l)=0$ if  $q>0$ or  $\lim X_t=-\i$ $\P$-a.s.), 
one has  $$|c_q(x,\l)|<+\i\qquad |b_q(x,\l)|<+\i$$ and  
$$A_q(x-,\l)b_q(x,\l)+\cc A_q(x,\l)c_q(x,\l)={1\over \psi_q(\l)}$$

For  $\Re(\l)<0$ (and for $\Re(\l)=0$ if $q>0$ or $\lim X_t=+\i$ $\P$-a.s.), 
one has  $$|\cc c_q(x,\l)|<+\i\qquad |\cc b_q(x,\l)|<+\i$$ and  
$$A_q(x-,\l)\cc d_q(x,\l)+\cc A_q(x,\l)\cc c_q(x,\l)={1\over \cc \psi_q(\l)}$$
\end{lemma}

{\bf Proof } Denote  $S_x:=\inf\{n, Z_n\in [-x,x[\}$, clearely $S_x=\tau_x+1$. 
One has  
$$\Q_q^{\uparrow}(\tau_x=0, Z_{S_x}\in dy)+\sum_{n=1}^{+\i} \Q_q^{\uparrow}( Z_1\in dx_1\dots Z_n\in dx_n \tau_x=n, Z_{S_x}\in dy)$$
$$=\Q_q^{\uparrow}(\tau_x=0, Z_{S_x}\in dy)+\sum_n 1_{x_n\not\in [-x,x[}1_{y\in [-x,x[} \Q_q^{\uparrow}( Z_1\in dx_1\dots Z_n\in dx_n, Z_{n+1}\in dy)$$
$$=H_q(dy)1_{y\in [0,x[}[1+\sum_{n  \hbox{ is even},n\geq 2}1_{ x_n\in ]-\i,-x[}\Q_q^{\uparrow}( Z_1\in dx_1\dots Z_n\in dx_n){1\over H_q((-x_n)^-)}]$$
$$+\cc H_q(-dy)1_{y\in [-x,0]}[\sum_{n  \hbox{ is odd }} 1_{ x_n\in [x,+\i[}\Q_q^{\uparrow}( Z_1\in dx_1\dots Z_n\in dx_n ){1\over \cc H_q(x_n)}]$$

Again, the first identity is a consequence of the fact that  the sequence  $|Z_n|$ is nonincreasing
and the second from the value of the transition kernel of the Markov chain given in lemma \ref {MC}.

Now, let us  check that  the sigma field ${\mathcal G}_{\tau_x}:=\sigma(Z_n1_{n\leq \tau_x}; n\geq 1)$ is $\sigma$-finite. Indeed the variable $Z_11_{\tau_x\not=0}+ 1_{\tau_x=0}$ is positive and is ${\mathcal G}_{\tau_x}$-measurable. The  variable  $Z_11_{\tau_x\not=0}$ has law $H_q(dy)1_{y\in [x,+\i[}$ which is sigma finite and the event $\{\tau_x=0\}$ has measure $H_q(x^-)$.
So one can deduce the conditionnal law of the variable $Z_{T_x}$ on  ${\mathcal G}_{\tau_x}$ from the preceeding identities 
$$\Q_q^{\uparrow}(Z_{S_x}\in dy\vert  {\mathcal G}_{\tau_x})=1_{\tau_x\hbox{ is odd}} {\cc H_q(-dy)\over \cc H_q(x)}1_{y\in [-x,0]}
 + 1_{\tau_x\hbox{ is even}} {H_q(dy)\over H_q(x^-)}1_{y\in [0,x[}$$
On the other hand, one has 
$$A_q(x^-,\l)=\Q_q^{\uparrow}(e^{-\l F};M<x)=\Q_q^{\uparrow}(e^{-\l \sum_1^{+\i} Z_n}; Z_1<x)=\int _{[0,x[}\Q_q^{\uparrow}(e^{-\l \sum_1^{+\i} Z_n}\vert Z_1=y)H_q(dy)$$
 and similarly, 
  $$\cc A_q(x,\l)=\int _{[-x,0]}\Q_q^{\downarrow}(e^{-\l \sum_1^{+\i} Y_n}\vert Y_1=y)\cc H_q(-dy)$$
Since  $(Z_n)$ is a Markov chain under  $\Q_q^{\uparrow}$ with same transition kernel as the chain  $(Y_n)$ under  $\Q_q^{\downarrow}$, and 
since  $S_x$ is a stopping time of the chain $(Z_n)$ and since the sigma-field  $\sigma( Z_n1_{n\leq S_x})$ contains  ${\mathcal G}_{\tau_x}$, one can deduce from previous identities that  :

$$\Q^{\uparrow}_q( e^{-\l \sum _{S_x}^{+\i} Z_n}\vert {\mathcal G}_{\tau_x})=1_{\tau_x\hbox{ is odd}}{\cc A_q(x,\l)\over \cc H_q(x)} +1_{\tau_x\hbox{ is even}}{A_q(x^-,\l)\over H_q(x^-)} $$
Finaly, one gets 
 $${1\over\psi_q(\l)}=\Q^{\uparrow}_q( e^{-\l F})= \Q^{\uparrow}_q( e^{-\l\sum_1^{+\i}Z_n})=
  \Q^{\uparrow}_q( e^{-\l\sum_1^{\tau_x}Z_n} e^{-\l\sum_{T_x}^{+\i}Z_n})$$
  $$=\Q_q^{\uparrow}( e^{-\l\sum_1^{\tau_x}Z_n}1_{\tau_x\hbox{ is even}}){A_q(x^-,\l)\over H_q(x^-)} +\Q_q^{\uparrow}( e^{-\l\sum_1^{\tau_x}Z_n}1_{\tau_x\hbox{ is odd}})  {\cc A_q(x,\l)\over \cc H_q(x)}$$
  $$= b_q(x,\l) A_q(x^-,\l)+c_q(x,\l)\cc A_q(x,\l)$$
  
  The single case for which one has to check that $c_q(x,\l)$ and  $d_q(x,\l)$ are actually finite is the case $q=0$ and $\Re(\l)>0$. It is enough  to check the property for positive real $\l$. In this case the real numbers $A_0(x^-,\l)$, $\cc A_0(x,\l)$ are positive and ${1\over\psi_0(\l)}$ is finite. Thus the previous identity allows us to conclude that the (positive but possibly infinite a priori) $c_0(x,\l)$ and $ b_0(x,\l)$ are actualy finite.

 One gets the second part of the lemma similarly.\qed

\begin{lemma} \label {Edc}
For  $\Re(\l)>0$ (and for  $\Re(\l)=0$ if  $q>0$ or if  $\lim X_t=-\i$ $\P$-a.s.),  one has
$$c_q(x,\l)=\int_{[x,+\i[} e^{-\l y} b_q(y,\l) {H_q(dy)\over \cc H_q(y)}$$
$$b_q(x,\l) =1+ \int_{]x,+\i[} e^{\l y} c_q(y,\l) {\cc H_q(dy)\over H_q(y^-)}$$
For $\Re(\l)<0$ (and for  $\Re(\l)=0$ if $q>0$ or $\lim X_t=+\i$ $\P$-a.s.), one has
$$\cc c_q(x,\l)=\int_{]x,+\i[} e^{\l y} \cc b_q(y,\l) {\cc H_q(dy)\over H_q(y^-)}$$
$$\cc b_q(x,\l) =1+ \int_{[x,+\i[} e^{-\l y} \cc c_q(y,\l) {H_q(dy)\over \cc H_q(y)}$$
\end{lemma}

{\bf Proof}
Extend the transition kernel of $(Z_n)$, $R_q$, by setting $R_q(y,\{-\i\}):=1-R_q(y, \R^*)$ for  $y\not =0$. One gets from the duality identity 
$1_{yz\not=0}U_q(dz)P_q(z,dy)=1_{yz\not=0}U_q(dy)R_q(y,dz)$
$$1_{y\in \R^*}U_q(y, \R^*)U_q(dy)=1_{y\in \R^*}\int_{ \R^*}U_q(dz)P_q(z,dy)=1_{y\in \R^*}\Q^{\uparrow}_q (\sum_2^{+\i}1_{Z_n\in dy})=1_{y\in \R^*}(U_q(dy) -H_q(dy))$$
and deduce :
$$R_q(y,\{-\i\})U_q(dy)=H_q(dy) \quad \hbox {on}\quad \R^*$$ 
For $y\in \R^*$, denote $p^y$ the distribution of the Markov chain starting from $y$ and with transition kernel $R_q$,  $(U_n)$ the canonical process  of the space $\left(\R^*\cup\{-\i\}\right)^{\N}$ where $-\i$ is the cemetery point and  $\xi$ the life time of the process $(U_n)$.
One can deduce from lemma \ref{MR}, the identities :

$$c_q(x,\l)={1\over \cc H_q(x)}\Q^{\uparrow}_q(e^{-\l\sum_{n=1}^{\tau_x} Z_n};\tau_x \hbox{ is odd})=
\int_{[x,+\i[}p^y( e^{-\l \sum_1^{\xi}U_n}) {U_q(dy) \over \cc H_q(y)}$$
and 
$$b_q(x,\l)={1\over H_q(x^-)}\Q^{\uparrow}(e^{-\l\sum_{n=1}^{\tau_x} Z_n};\tau_x \hbox{ is even})$$
$$={1\over H_q(x^-)}[\Q^{\uparrow}_q({\tau_x}=0)+\int_{]-\i,-x]}p^y( e^{-\l \sum_1^{\xi}U_n}). U_q(dy) {H_q(x^-)\over H_q((-y)^-)}]$$
$$=1+ \int_{]-\i,-x[}p^y( e^{-\l \sum_1^{\xi}U_n}). {U_q(dy) \over H_q((-y)^-)}$$ 

In the other hand,  one gets (measures involved here are on $\R^*$) 
$$p^y( e^{-\l \sum_1^{\xi}U_n})U_q(dy) =e^{-\l y} [\int _{ [-\i,+\i[\backslash \{0\}}p^z( e^{-\l \sum_1^{\xi}U_n})R_q(y,dz)]U_q(dy)$$  
$$=e^{-\l y}[R_q(y,\{-\i\})U_q(dy) +\int _{ ]-\i,+\i[\backslash \{0\}}p^z( e^{-\l \sum_1^{\xi}U_n})U_q(dz)P_q(z,dy)]$$
$$=e^{-\l y}[1+\int _{]-\i,-y[}p^z( e^{-\l \sum_1^{\xi}U_n}){U_q(dz)\over H_q((-z)^-)}] H_q(dy)=e^{-\l y} b_q(y,\l)H_q(dy)$$
The first identity comes from the Markov property of $(U_n)$ under $p^y$ with transition kernel $R_q$, the second one from the duality identity 
 $R_q(y,dz)U_q(dy)=U_q(dz)P_q(z,dy)$, the
third one from the value of 
 $P_q(z,dy)$ given in lemma \ref {MC} and the identity  $R_q(y,\{-\i\})U_q(dy)=H_q(dy)$, and the fourth one from the expression of  $b_q(y,\l)$ given below. 
Replacing the preceeding identity into the expression of   $c_q(x,\l)$ given previously, one gets 

$$c_q(x,\l)=\int_{]x,+\i[} e^{-\l y} b_q(y,\l) {H_q(dy)\over \cc H_q(y)}$$

One proves the other assertions of the lemma in a similar way. \qed

\begin{lemma}  \label {caC} For all $x\in ]0,+\i[$, $\Re(\l)> 0$ (and $\Re(\l)=0$ if $q>0$ or $\lim_{t\to +\i} X_t=-\i$, $\P$-ps) one has
$$\psi_q(\l)c_q(x,\l) =C_q(x^-,\l)\quad\hbox{and}\quad \psi_q(\l)b_q(x,\l)=B_q(x,\l)$$
For all $x\in ]0,+\i[$, $\Re(\l)< 0$ (and $\Re(\l)=0$ if $q>0$ or $\lim_{t\to +\i} X_t=+\i$, $\P$-ps), one has
$$\cc \psi_q(\l)\cc c_q(x,\l) =\cc C_q(x,\l)\quad\hbox{and}\quad \cc \psi_q(\l)\cc b_q(x,\l)=\cc B_q(x^-,\l)$$
\end{lemma} 
{\bf Proof} 
Remember that the time $U_x$ is defined as 
$$U^x=\inf\{t;  S_t-I_t\geq x \hbox { and  }S_t=X_t \}\wedge \inf\{t;  S_t-I_t>x \hbox { and }X_t=I_t \}$$
Let us denote $B^x$ the event $B^x= \{X_{U^x}=S_{U^x}, U^x<+\i\}$
and denote  $\hat X$ the process  $\hat X_t(w)=F-X_{t-}(w)$. Add a $\hat{} $ for the corresponding objects.
Denote 
 $X^{\uparrow}$ and  $X^{\downarrow}$ the processes $(X_{t+\rho}-m; t\geq 0\}$ and  $(X_{t+\sigma}- M;  t\geq 0\}$ respectively.

A quick look at the picture will convince the reader that the event $\{\tau_x(X^{\uparrow})\hbox{ is odd}\}$ is equal to the event $\hat B^x$ and on this event, one has :
$\sum_1^{\tau_x(X^{\uparrow})} Z_n(X^{\uparrow})= \hat M-\hat I_{\hat U_x}$.
Take any  event $A$ such that $0<\Q^{\downarrow}_q(A)<+\i$ and get 
$$\Q_q^{\downarrow}(A)\cc H_q(x)c_q(x,\l)=\Q_q^{\downarrow}(A)\Q_q^{\uparrow}(e^{-\l\sum_{n=1}^{\tau_x} Z_n};\tau_x \hbox{is odd})$$
$$=\Q_q(A(m-X_{(\rho-t)^-}; t\geq 0); e^{-\l\sum_{n=1}^{\tau_x(X^{\uparrow})} Z_n};\tau_x (X^{\uparrow})\hbox{ is odd})$$
$$=\Q_q( A(\hat X^{\downarrow}),e^{-\l (\hat M-\hat I_{\hat U_x})}; \hat B_x)=\Q_q(e^{-\l (M-I_{U_x})}; B_x, A(X^{\downarrow}))$$ 
The first identity follows from the definition of $c_q$,  the second one from proposition \ref{DT}, the 
third one from what we have just  said, 
and  the fourth one from the identity in law of $\hat X$ and  $X$ under $\Q_q$.

On the other hand, $U_x$ is a stopping time and it is smaller than  $\sigma$ on the event $B^x$, therefore we get 
 (denote  $\th^{\circ}_{U_x}$ the path $(X_{t+U_x}-X_{U_x}; t\geq 0)$)
 $$M-I_{U^x}= (X_{U_x}-I_{U_x})+M\circ \th^{\circ}_{U_x} \quad \hbox{and } \quad A(X^{\downarrow})=A(X^{\downarrow}\circ \th^{\circ}_{U_x} )\hbox{ on  } B^x$$
We deduce then the following identities from the property of independent increments at
the stopping
 time $U_x$ :

$$\Q_q( e^{-\l (M-I_{U_x})}; B^x; A(X^{\downarrow}))=\Q_q( e^{-\l (X_{U_x}-I_{U_x})}; B^x; e^{-\l M\circ \th^{\circ}_{U_x} }; A(X^{\downarrow}\circ \th^{\circ}_{U_x} ))$$
$$=\P_q( e^{-\l (X_{U_x}-I_{U_x})}; B^x)\Q_q(e^{-\l M};A(X^{\downarrow}))$$

First identity of  proposition \ref {A} for $\l=0$, $\mu_1=\l$ and $q_1=q_2=q$ 
gives us  $$\P_q( e^{-\l (X_{U_x}-I_{U_x})}; B^x)=\cc H_q(x)C_q(x^-,\l)$$
and proposition \ref{DT} gives us :
$$\Q_q(e^{-\l M};A(X^{\downarrow}))=\Q_q^{\uparrow}(e^{-\l F}) \Q_q^{\downarrow}(A)={1\over\psi_q(\l)}\Q_q^{\downarrow}(A)$$

Putting together the previous identities and simplifying by  $\Q_q^{\downarrow}(A)\cc H_q(x)$ on gets the following 
$$C_q(x,\l)= \psi_q(\l)c_q(x,\l)$$

In order to prove the second identity of the lemma, let us remind that
$$T^s_x=\inf\{ t; X_t-S_t<- x\}\qquad {1\over B_q(x,\l)}=\P_q(\int_{[0,{T^s_x}[}e^{-\l S_t} L^s(dt))\qquad A_q(x^-,\l)=\P_q(\int_{[0,U^x[} e^{-\l S_t} L^s(dt))$$
 
 Notice  that if $U^x<+\i$ and $I_{U_x}=X_{U_x}$  then  $T^s_x=U^x$ 
and if  $U^x=+\i$ then  $T^s_x=+\i$.
In other cases, one has  $X_{U^x}=S_{U^x}$  and (still denoting  $\th^{\circ}_{U_x}$ the path  $(X_{s+U^x}-X_{U^x}; s\geq 0)$)  :
$$T^s_x=U^x+T^s_x\circ \th^{\circ}_{U_x}\qquad  L^s(w,dt+U^x)1_{t>0}=L^s(\th^{\circ}_{U_x}(w),dt)$$
$$ S_{U^x+t}=X_{U^x}+S_t\circ \th^{\circ}_{U_x}\quad \hbox{ for all } t\in [0,+\i[$$
One gets the identity of random variables : 
$$\int_{[0,{T^s_x}[}e^{-\l S_t}L^s(dt)=\int_{[0,U^x[}e^{-\l S_t}L^s(dt) +1_{\{X_{U^x}=S_{U^x};U^x<+\i\}}e^{-\l X_{U^x}}[\int_{[0,{T^s_x}[}e^{-\l S_t}L^s(dt)]\circ \th^{\circ}_{U_x}$$
Take the expectation with respect to $\P_q$  and use the property of independence of increments at time $U_x$ to get  

$${1\over B_q(x,\l)}=A_q(x^-,\l)+ \P_q(e^{-\l X_{U^x}}; X_{U^x}=S_{U^x};U^x<+\i){1\over B_q(x,\l)}$$

In the other hand, use first identity of proposition \ref{A} for $\mu_1=\l$ and $q_1=q_2=q$  and get 
$$\P_q(e^{-\l X_{U^x}}; X_{U^x}=S_{U^x};U^x<+\i)=\cc A_q(x,\l)C_q(x^-,\l)$$
One deduces 
$$1=B_q(x,\l)A_q(x^-,\l)+\cc A_q(x,\l)C_q(x^-,\l)$$
When comparing with the identity obtained at lemma \ref{bc} :
$A_q(x^-,\l) b_q(x,\l)+\cc A_q(x,\l)c_q(x,\l)={1\over \psi_q(x,\l)}$ 
and the identity already obtained $C_q(x^-,\l)=\psi_q(x,\l)c_q(x,\l)$, 
we get  the following 
$$B_q(x,\l)=\psi_q(\l) b_q(x,\l)$$

A similar proof works to obtain the two other identities of the lemma.\qed

Properties 2 and  3 of theorem \ref {ED} follow from  lemmas \ref {bc}, \ref{Edc}  and \ref {caC}.

\section{Proof of theorem \ref{RH}}
\subsection {Two lemmas}
\begin{lemma}\label {ST} For all $x\in ]0,+\i[$ and $iu\in i\R$, one has 
$$B_q(x,iu)\cc B_q(x^-,iu) - C_q(x^-,iu)\cc C_q(x,iu)=\phi_q(iu)$$
$$\cc B_q(x,iu)= (\phi(iu)+q) A_q(x,iu)+ C_q(x,iu)$$
$$B_q(x,iu)=(\phi(iu)+q)\cc A_q (x,iu)+\cc C_q(x,iu)$$

\end{lemma}
{\bf Comments} According to  theorem \ref{ED}  the three pairs of functions $(A_q(x^-,iu),\cc A_q(x,iu))$, $(C_q(x^-,iu),B_q(x,iu))$ and $(\cc B_q(x^-,iu),\cc C_q(x,iu))$ satisfy the same differential equation in $x$,  in the sense of distribution theory. If the coefficients of this equation were  suffisiently regular, the second and the third identities of  lemma 6.1 would be an easy  consequence of 2- dimensionality of the space of solution of this
 differential  equation and the first identity would come from the wronskian identity. In the next proof, we check that these results still hold in our setting by using Stieljes integral calculus.

{\bf Proof} First, notice that every complex valued function defined on $[0,+\i]$ which is bounded and have bounded variations 
on every half line $[x,+\i]$  (respectively on every interval $[0,x[ $) where  $x$ is a positive real can be written  $\int _{]x,+\i]} f(t)\nu(dt)$ (resp. $\int _{[0,x]} f(t)\nu(dt)$) if it is right continuous or $\int _{[x,+\i]} f(t)\nu(dt)$(resp. $\int _{[0,x[} f(t)\nu(dt)$) if it is left continuous  where $\nu$ is a positive measure on $[0,+\i]$ and $f$ is a complex valued function. We shall denote $\mu]x,+\i]$ or $\mu[x,+\i[$(resp. $\mu[0,x]$ or $\mu[0,x[$) these functions in the sequel and $\mu(dx)$ will be the associated complex measure ($\mu(dx)=f(x)\nu(dx)$).
Here are few easy facts about Stieljes integral calculus.
 
 An easy application of Fubini's theorem gives us, for all positive  $x$,

$$\mu_1[x,+\i]\mu_2]x,+\i]=\int_{]x,+\i]}\mu_1[y,+\i]\mu_2(dy)+\int_{[x,+\i[}\mu_2]y,+\i]\mu_1(dy)\eqno{(1)}$$

If the complex measure $\mu_1(dx)$ is also integrable in the neighborhood of $0$, one has :
$$\mu_1[0,x[\mu_2]x,+\i]=\int_{]x,+\i]}\mu_1[0,y[\mu_2(dy)-\int_{[x,+\i[}\mu_2]y,+\i]\mu_1(dy)$$
[To  establish  this one, add it with previous one].
When  adding the term $\sum_{y\in ]x,+\i[}\mu_1(\{y\})\mu_2(\{y\})$ in both integrals, one gets :
$$\mu_1[0,x[\mu_2]x,+\i]=\int_{]x,+\i]}\mu_1[0,y]\mu_2(dy)-\int_{[x,+\i[}\mu_2[y,+\i]\mu_1(dy)+\mu_2(\{x\})\mu_1(\{x\})$$
And by regularisation on the right, one has 
$$\mu_1[0,x]\mu_2]x,+\i]=\int_{]x,+\i]}\mu_1[0,y]\mu_2(dy)-\int_{]x,+\i[}\mu_2[y,+\i]\mu_1(dy)\eqno{(2)}$$
Applying the previous equation to $\mu_1[0,x]:={1\over \mu_2]x,+\i]}$ when ${1\over \mu_2]x,+\i]}$ has bounded variations on $[0,x]$ for every  $x$, one gets the identity of complex measures 
$$\mu_1[0,y]\mu_2(dy)=\mu_2[y,+\i]\mu_1(dy)\quad \hbox{ on }\quad ]0,+\i[$$
and one deduces $$\mu_1[0,x]={1\over \mu_2]x,+\i]} \Longrightarrow \mu_1(dy)={\mu_2(dy)\over \mu_2]y,+\i]\mu_2[y,+\i]}\quad \hbox{ on }\quad ]0,+\i[\eqno{(3)}$$

Let us establish now the  first identity of the  lemma \ref {ST}.
When applying identity $(1)$ to 
$$\mu_1[x,+\i]:=\cc B_q(x^-,iu)=\cc \psi_q(iu)+ \int_{[x,+\i[} \cc C_q(y,iu) e^{-iuy} {H_q(dy)\over \cc H_q(y)}$$
and 
$$\mu_2]x,+\i]:=B_q(x,iu)=\psi_q(iu)+ \int_{]x,+\i[} C_q(y^-,iu) e^{iuy} {\cc H_q(dy)\over H_q(y^-)}$$
one gets 
$$ \cc B_q(x^-,iu)B_q(x,iu)=\cc \psi_q(iu) \psi_q(iu)+ \int_{]x,+\i[}\cc B(y^-,iu).C(y^-,iu) e^{iuy} {\cc H_q(dy)\over H_q(y^-)}$$
$$+ \int_{[x,+\i[}B(y,iu).
\cc C(y,iu) e^{-iuy} {H_q(dy)\over \cc H_q(y)}$$

when applying this identity to 
$$\mu_1[x,+\i]:=C_q(x^-,iu)= \int_{[x,+\i[} B_q(y,iu) e^{-iuy} {H_q(dy)\over \cc H_q(y)}$$
and 
$$\mu_2]x,+\i]:=\cc C_q(x,iu)= \int_{]x,+\i[} \cc B_q(y^-,iu) e^{iuy} {\cc H_q(dy)\over H_q(y^-)}$$

one gets 
 $$C_q(x^-,iu)\cc C_q(x,iu)= \int_{]x,+\i[}C_q(y^-,iu).\cc B_q(y^-,iu) e^{iuy} {\cc H_q(dy)\over H_q(y^-)}+ \int_{[x,+\i[}\cc C_q(y,iu).B_q(y,iu) e^{-iuy} {H_q(dy)\over \cc H_q(y)}$$
 When substracting, these equations, one obtains 
 $$ \cc B_q(x^-,iu)B_q(x,iu)-C_q(x^-,iu)\cc C_q(x,iu)=\cc \psi_q(iu) \psi_q(iu)=\phi(iu)+q$$
 
Let us now establish the second identity of lemma \ref{ST}.
When multiplying by 
 $\phi(iu)+q$ the identity obtains in theorem \ref{ED} 
 ($\cc A_q(x,iu)\cc B_q(x^-,iu)+A_q(x^-,iu)\cc C_q(x,iu) =1$) and substracting it from the  previous one, one gets  
$$ \cc B_q(x^-,iu)[B_q(x,iu)-(\phi(iu)+q)\cc A_q(x,iu)]- \cc C_q(x,iu) [C_q(x^-,iu)+(\phi(iu)+q)A_q(x^-,iu)]=0$$
Multiplying this last equation by the measure ${1\over \cc B_q(x^-,iu) \cc B_q(x,iu)}e^{-iux} {H_q(dx)\over \cc H_q(x)}$, one obtains 

$$ \Bigl[{1\over \cc B_q(x,iu)}\Bigr]\Bigl [(B_q(x,iu)-(\phi(iu)+q)\cc A_q(x,iu))e^{-iux} {H_q(dx)\over \cc H_q(x)}\Bigr]$$
$$-\Bigl[C_q(x^-,iu)+(\phi(iu)+q)A_q(x^-,iu)\Bigr]\Bigl[{ \cc C_q(x,iu)\over  \cc B_q(x^-,iu) \cc B_q(x,iu)}e^{-iux} {H_q(dx)\over \cc H_q(x)}\Bigr]=0\eqno{(4)}$$

Put $q>0$ and let us justify that equation (4) is of the form  $\mu_1[0,x]\mu_2(dx)-\mu_2[x,+\i]\mu_1(dx)=0$ where 
$$\mu_1[0,x]:={1\over \cc B_q(x,iu)}\quad\hbox{et}\quad \mu_2[x,+\i]:=C_q(x^-,iu)+(\phi(iu)+q)A_q(x^-,iu)$$
Indeed, bounded variations on $[x,+\i]$ of $C_q(x^-,iu)+(\phi(iu)+q)A_q(x^-,iu)$ is justified as follows :
The measures $\Q^{\uparrow}_q(dw)$ and $H_q (dx)$ are finite because $q$ is positive, thus the function  $A_q(x^-,iu)=\Q_q^{\uparrow}(M<x;e^{-iu F})$ tends to 
$\Q_q^{\uparrow}(e^{-iu F})={1\over\psi_q(iu)}$ when $x$ goes to $+\i$. In the other hand, the function of $x$, ${\cc A_q(x,iu)\over \cc H_q(x)}$ is bounded by $1$ and so it is $H_q(dx)$ integrable. We deduce from these facts and from first  identity of theorem \ref{ED} ($A_q(x^-,iu)=H_q(0)+\int_{]0,x[}e^{-iuy} {\cc A_q(y,iu)\over \cc H_q(y)} H_q(dy)$), the next one 
$$A_q(x^-,iu)={1\over\psi_q(iu)}-\int_{[x,+\i[} e^{-iuy}\cc A_q(y,iu){H_q(dy)\over \cc H_q(y)}$$
Consequently, this equation and the second identity of theorem \ref {ED} ($C_q(x^-,iu)=\int_{[x,+\i[}e^{-iuy}  B_q(x,iu){ H_q(dy)\over \cc H_q(y)}$)
gives us that  
$$\mu_2[x,+\i]:=C_q(x^-,iu)+(\phi(iu)+q)A_q(x^-,iu)$$
has bounded variations on $[x,+\i]$ for every positive real $x$
and $$\mu_2(dx)=(B_q(x,iu)- (\phi(iu)+q)\cc A_q(x,iu))e^{-iux} {H_q(dx)\over \cc H_q(x)} \quad \hbox{on}\quad ]0,+\i[$$
In the other hand, bounded variations on $[0,x]$ of ${1\over \cc B_q(x,iu)}$ follow from the identity 
$${1\over \cc B_q(x,iu)}=\P(\int_{[0,T_i^x[}e^{-iu I_t-qt}L^i(dt))$$
Apply formula $(3)$ to $\mu_1[0,x]:={1\over B_q(x,iu)}={1\over \mu_2 ]x,+\i]}$ and use the identity of theorem \ref {ED}
($\cc B_q(x^-,iu)=\int_{[x,+\i[}e^{-iuy} \cc  C_q(x,iu){ H_q(dy)\over \cc H_q(y)}$) and get 
 $$\mu_1(dx)= {\cc C_q(x,iu) \over  \cc B_q(x^-,iu) \cc B_q(x,iu)}e^{-iux}{H_q(dx)\over \cc H_q(x)}\quad \hbox{on}\quad ]0,+\i[$$
 Now, one can  integrate equation $(4)$ over  $]x,+\i]$ and use equation $(2)$ and  get that the  function 
$$\mu_1[0,x]\mu_2]x,+\i]={C_q(x,iu)+(\phi(iu)+q)A_q(x,iu)\over \cc B_q(x,iu)}$$ is a constant function of $x$ on $]0,+\i[$. 
0n the other hand, 
$$\lim_{x\to +\i}C_q(x,iu)=0\qquad \lim_{x\to +\i}A_q(x,iu)={1\over \psi_q(iu)}\qquad \lim_{x\to +\i} \cc B_q(x,iu)=\cc \psi_q(iu)$$ Thus 
$$\lim_{x\to +\i}{C_q(x,iu)+(\phi(iu)+q)A_q(x,iu)\over \cc B_q(x,iu)}={\phi(iu)+q\over \psi_q(iu)\cc \psi_q(iu)}=1$$
This enables us to deduce the second identity of the lemma \ref {ST}  $$\cc B_q(x,iu)= C_q(x,iu)+(\phi(iu)+q)A_q(x,iu)$$

The third identity can be proved similarly  or by regularizing the second one at the left  and integrating it over $]x,+\i]$ with respect to the measure 
$e^{iux}{\cc H_q(dx)\over H_q(x^-)}$.
 To get the second and the third identity when $q=0$, take a limit.\qed

Enlarge the probability space $(\Omega,{\mathcal F},\P)$ in order that it contains a random variable $\xi_q$ independent of the canonical process  $X$ and which has  an exponential distribution with parameter $q$ and still denote by the same notation this enlarged space. 
For $t>0$, denote 
$$D^i_t=\inf \left([t,+\i[\cap \hbox{supp }( L^i)\right)$$

\begin{lemma}\label  {LM} For  $\Re(\l)\geq 0$, one gets  
$$1-\cc \psi_0 (\l)\cc A_q(x,\l)=\P(e^{-\l I_{D^i_{U^x\wedge \xi_q}}};D^i_{U^x\wedge \xi_q}<+\i)$$

\end{lemma}
{\bf Proof}. Denote $\th^{\circ}_{D^i_{U^x\wedge \xi_q}}$ the path  $(X_{t+D^i_{U^x\wedge \xi_q}}-X_{D^i_{U^x\wedge \xi_q}}; t\geq 0)$, one
has 
$$\int_{[0,\i[} e^{-\l I_t} L^i(dt)= \int _{[0,U^x\wedge \xi_q[} e^{-\l I_t} L^i(dt)+1_{D^i_{U^x\wedge \xi_q}<+\i}. e^{-\l I_{D^i_{U^x\wedge \xi_q}}}[\int_{[0,+\i[} e^{-\l I_t}L^i(dt)]\circ \th^{\circ}_{D^i_{U^x\wedge \xi_q}}$$
Take the expectation of this equation with respect to $\P$, use the independent increments at stopping time $D^i_{U^x\wedge \xi_q}$ and get 

$${1\over\cc \psi_0(\l)} = \cc A_q(x,\l)+ \P(e^{-\l I_{D^i_{U^x\wedge \xi_q}}};D^i_{U^x\wedge \xi_q}<+\i ){1\over\cc \psi_0(\l)}$$
Thus  
$$1- \cc\psi_0(\l) \cc A_q(x,\l)=  \P(e^{-\l I_{D^i_{U^x\wedge \xi_q}}};D^i_{U^x\wedge \xi_q}<+\i )$$

\qed

\subsection{The matrix $M_q$ satisfies properties of theorem \ref {RH}}

Property 1)  is obvious. Property  2) follows from  the second and the third identities of lemma \ref {ST}. Let us check property 3) :

The functions $\l\mapsto A_q(x^-,\l)$,  $\l\mapsto \cc A_q(x,\l)$, $\l\mapsto C_q(x^-,\l)$, $\l\mapsto \cc C_q(x,\l)$ are respectively Laplace transforms of finite measures supported respectively by sets $[0,x[$, $[-x,0]$, $[x,+\i[$ and $]-\i,x]$. Consequently, 
the functions $\l\mapsto A_q(x^-,\l)$, $\l\mapsto e^{-\l x} \cc A_q(x,\l)$, $\l\mapsto e^{\l x}C_q(x^-,\l)$ are bounded on the half plane  $\{\Re(\l)>0\}$ and the functions $\l\mapsto e^{\l x} A_q(x^-,\l)$, $\l\mapsto \cc A_q(x,\l)$, $\l\mapsto e^{-\l x}\cc C_q(x,\l)$ are bounded on the half plane  $\{\Re(\l)<0\}$.

 The function $\l \mapsto B_q(x,\l)$ (resp. $\l \mapsto B_q(x^-,\l)$) is the L\'evy exponent of a subordinator (resp. opposite of a subordinator) 
killed at an exponential independent time, thus the function $\l \to {B_q(x,\l)\over |\l| +1}$ (resp.  $\l \to {\cc B_q(x^-,\l)\over |\l| +1}$) is bounded on the half plane $\{\Re(\l)<0\}$ (resp. $\{\Re(\l)>0\}$).

Let us check now property  4).  The function $\l\mapsto e^{\l x} A_q(x^-,\l)$ is the Laplace transform of a finite mesure supported by $[-x,0[$ thus  
$$\lim_{\Re(\l)\to -\i}e^{\l x} A_q(x^-,\l)=0$$

Note that the random variable $I_{D^i_{U^x\wedge \xi_q}}$ is negative (If $\P$ is not the distribution of a compound Poisson process,  $I_t$ is negative on $\hbox{supp }L^i\backslash \{0\}$ as we have already said in part I of the paper  and if $\P$ is the distribution of a compound Poisson process, it is still the case by definition of the local time $L^i$), thus the limit $$\lim_{\Re(\l)\to -\i}\cc\psi_0(\l)\cc A_q(x,\l) =1$$  follows from lemma \ref {LM}.

The limit $$\lim_{\Re(\l)\to -\i}{ \cc B_q(x^-,\l)\over \cc \psi_0(\l)}=1$$ follows from the two previous limits, the boundedness of $\l\mapsto e^{-\l x}\cc C_q(x,\l)$ and  the identity $\cc A_q(x,\l)\cc B_q(x^-,\l)+A_q(x^-,\l)\cc C_q(x,\l)=1$ stated in theorem \ref {ED}.

If  $]-\i,0[$ is irregular then for every positive stopping time $T$,  either $\P$ is the distribution of compound Poisson  process and 
the set  $\{s; \inf _{0\leq u\leq s}X_{u+T}-X_{T^-}=X_{s+T}- X_{T^-}\}$ is ($\P$-a.s.) a discrete union of intervals where the process $X$ is constant, either  the set  $\{s; \inf _{0\leq u\leq s}X_{u+T}-X_{T^-}=X_{s+T}-X_{T^-}\}$ is discrete.  
One can easily deduce that the set  $\{t; I_t=X_t\}$ has the same property $\cc N$-a.s. 
On the other hand, the  time  $T_x$  is defined to be the first time when $X_t$ is strictly smaller than  $-x$. One  deduces that $X_{T_x} <-x$ on $T_x<+\i$ $\cc N$-as, and  
$\l\mapsto \cc C_q(x,\l)=\cc N( e^{-\l X_{T_x}-qT_x};T_x<+\i)$ is the Laplace transform of a finite measure supported by $]-\i,x[$, thus 

$$\lim_{\Re(\l)\to -\i}e^{-\l x}\cc C_q(x,\l)=0$$

\subsection{Proof of the uniquiness part of theorem \ref{RH} } 
Note first that the
identities $$A_q(x^-,\l)B_q(x,\l)+\cc A_q(x,\l) C_q(x^-,\l)=1$$
 and 
 $$\cc A_q(x,\l)\cc B_q(x^-,\l)+A_q(x^-,\l)\cc C_q(x,\l)=1$$ stated in theorem \ref {ED} imply that $\det M_q(x,\l)=1$ for every $\l \in \C\backslash i\R$.  
Let now $N(\l)$ be another matrix having the same property as $M_q(\l)$ and let us  check  that  $N(\l)[M_q(\l)]^{-1}$ is equal to the identity matrix.
For simplicity we omit the parameters $x$ and $\l$  in the notation  in the sequel. 

By property 1)  the matrix $N(\l)[M(\l)]^{-1}$ is analytic  on the two half-plane  $\{\Re(\l)>0\}$and $\{\Re(\l)<0\}$, and by property 2) it can be 
extended by continuity to every point of the imaginary axis. 

Property 3) allows us to state that $N(\l)M(\l)^{-1}$is bounded on every compact set of $C$.  Thus  a standard argument gives us that the extended matrix   is entire. 

Property 3) gives us also that the matrix 
${1\over (|\l|+1)^2} e^{\l {x\over 2} J} N(\l)[M(\l)]^{-1} e^{-\l {x\over 2} J}$ is bounded on $\C$ (remember that $J$ denotes the matrix  $\left( \begin{array} {clcr} 1& 0\\ 0 & -1\end{array} \right)$), 
thus the components of matrix $e^{\l {x\over 2} J} N(\l)[M(\l)]^{-1} e^{-\l {x\over 2} J}$  are polynomials. 

In the other hand, if $]-\i,0[$ is irregular then property 4) gives us $e^{\l {x\over 2} J} N(\l)[M(\l)]^{-1} e^{-\l {x\over 2} J}$ is equivalent, for  $\Re(\l)\to -\i$, to 
$$ \left( \begin{array} {clcr} \cc \psi_0(\l)& 0\\ 0 & {1\over \cc \psi_0(\l)}\end{array} \right)\left( \begin{array} {clcr} \cc \psi_0(\l)& 0\\ 0 & {1\over \cc \psi_0(\l)}\end{array} \right)^{-1}$$
One deduces that $e^{\l {x\over 2} J} N(\l)[M(\l)]^{-1} e^{-\l {x\over 2} J}$ is constantly equal to the identity matrix and $N(\l)=M(\l)$.

If $]-\i,0[$ is regular then the matrix 
$e^{\l {x\over 2} J} N(\l)[M(\l)]^{-1} e^{-\l {x\over 2} J}$ is equivalent, for  $\Re(\l)\to -\i$,  to 
$$ \left( \begin{array} {clcr} \cc \psi_0(\l)& 0\\ e^{-\l x}N_{21}(\l) & {1\over \cc \psi_0(\l)}\end{array} \right)\left ( \begin{array} {clcr} \cc \psi_0(\l)& 0\\ e^{-\l x}M_{21}(\l)  & {1\over \cc \psi_0(\l)}\end{array} \right)^{-1}= \left( \begin{array} {clcr}1& 0\\  {1\over \cc \psi_0(\l)}(e^{-\l x} N_{21}(\l)-e^{-\l x} M_{21}(\l)) & 1\end{array} \right)$$
The function $\cc \psi_0(\l)$ goes to $+\i$  because $\cc \psi_0(\l)$ is  the  L\'evy exponent of the opposite of a subordinator  (possibly 
with a finite life time) which is not a  compound Poisson process. On the other hand,  the function $e^{-\l x} N_{21}(\l)-e^{-\l x} M_{21}(\l)$ stays bounded. Thus $e^{\l {x\over 2} J} N(\l)[M(\l)]^{-1} e^{-\l {x\over 2} J}$ converges to the identity matrix again and $N(\l)=M(\l)$.

\section{Probabilistic interpretation of the identities between the six functions}
First let us introduce some terminology on Wiener-Hopf factorization.
Let us denote  
$(W,{\mathcal G},\Q)$ any probabilistic space. Let $Z$ be a real process defined on  $(W,{\mathcal G})$ with a life time that may be finite.
 If the distribution of $Z$ under $\Q$ is $\P_q$, we shall say that  $\phi(iu)+q$ is the L\'evy exponent of  $Z$ under  $\Q$ and its spatial Wiener-Hopf factors are the  functions $\cc \psi_q$ and  $\psi_q$.

Let  $ S=(0,S_1,\dots, S_n,\dots)$  be a sequence of real random variables defined on $(W,{\mathcal G})$ and having under $\Q$ the distribution of a random walk possibly   killed at an independent geometric time. 
We shall identifie its distribution to 
the distribution of the compound Poisson process with L\'evy exponent $1-\Q(e^{-iu S_1}; 1<\xi)$ ($\xi$ denotes the life time of $S$). We shall say that this function is the L\'evy exponent of the random walk $S$ and that  its spatial Wiener-Hopf factors are the ones of that compound Poisson process. One can easily check that these factors are respectively $1-\Q(e^{-iu S_{V_0}};V_0<\xi)$ and 
$1-\Q(e^{-iu S_{V^{0-}}};V^{0-}<\xi)$.  ($V_0:=\inf\{n\geq 1; S_n<0\} $) and  $V^{0-}:=\inf\{n\geq 1; S_n\geq 0\})$.

Let us now go back to our space $(\Omega,{\mathcal F},\P)$ and denote $\th^{\circ}_t$ the path $(X_{s+t}-X_t; s\geq 0)$ for every time $t$. The reader can easily convince him/herself of the  next assertions :

\subsection{Interpretation of the identity  $\cc B_q(x,iu)= (\phi(iu)+q)A_q(x,iu)+ C_q(x,iu)$}\
\par
 If  $]-\i,0[$ is regular then define  :

$$T_0:=0\qquad T_1:=T^x_i\qquad T_{n+1}=T_n +T^x_i\circ \th^{\circ}_{T_n} \qquad (T_{n+1}:=+\i \hbox{ if }T_n=+\i)$$
$$\tilde L^i(dt):=  \sum_{n=0}^{+\i} 1_{T_n\leq t<T_{n+1}}L^i( dt -T_n, \th^{\circ}_{T_n})$$
$$Y_t:=X_{\tilde L^{i,-1}_t}$$
where $\tilde L^{i,-1}$ denotes the right inverse of the continuous increasing function $t\mapsto \tilde L^i[0,t]$.

Under $\P_q$, $Y$ is a L\'evy process, its L\'evy exponent is  
 $\cc B_q(x,iu)-C_q(x,iu)$ and its Wiener-Hopf spatial factors are  $\cc \psi_q(iu)$ and  $\psi_q(iu) A_q(x,iu)$.

If  $]-\i,0[$ is irregular then let us define the sequence of stopping time  $(T_n; n\geq 1)$ as follows
$$T_1:=T_0^x=\inf \{t, X_t\not\in [0,x]\}\qquad T_{n+1}:=T_n+T^x_0\circ \th^{\circ}_{T_n}\qquad (T_{n+1}:=+\i \hbox{ if }T_n=+\i)$$

Under  $\P_q$, the sequence  $(0,X_{T_1},\dots, X_{T_n},\dots)$ is a random walk. Its L\'evy exponent is  $\cc B_q(x,iu)-C_q(x,iu)$ and its Wiener-Hopf spatial factors are  $\cc \psi_q(iu)$ and $\psi_q(iu) A_q(x,iu)$.

The interpretation of the identity  $B_q(x,iu)= (\phi(iu)+q)\cc A_q(x,iu)+ \cc C_q(x,iu)$ is similar.

\subsection{Interpretation of the identity   $(\phi(iu)+q)\cc A_q(x,iu)A_q(x^-,iu)+ C_q(x^-,iu)\cc A_q(x,iu)+ \cc C_q(x,iu)A_q(x^-,iu)=1$}\ \par 
This identity is obtained from 
$$B_q(x,iu)A_q(x^-,iu) +C_q(x^-,iu)\cc A_q(x,iu)=1\quad \hbox{ and } B_q(x,iu)= (\phi (iu)+q) \cc A_q(x,iu)+\cc C_q(x,iu)$$
Denote $(U_n; n\geq 1)$ the sequence of stopping times :
$$U_1:=U^x\qquad U_{n+1}:=U_n+U^x\circ\th^{\circ}_{U_n}, \quad (U_{n+1}:=+\i \hbox{ if } U_n=+\i) $$

Under  $\P_q$, the sequence  $(0, X_{U_1}, \dots, X_{U_n},\dots)$ is a random walk, its L\'evy exponent is   $1- C_q(x^-,iu)\cc A_q(x,iu)- \cc C_q(x,iu)A_q(x^-,iu)$ and its Wiener-Hopf spatial  factors are $\cc \psi_q(iu) \cc A_q(x^-,iu)$ and $\psi_q(iu)A_q(x^-,iu)$.

\subsection{Interpretation of the identity $\cc B_q(x^-,iu) B_q(x,iu)- C_q(x^-,iu)\cc C_q(x,iu)=\phi(iu)+q$}\ \par 

Let $S=(0,S_1,\dots,S_n,\dots)$ be a real valued random walk (with infinite life time) on any probability space $(W,{\mathcal G},\Q)$. 
Denote for all $x\in ]0,+\i[$, 
$$V_0^{x-}=\inf\{n\geq 1; S_n\not\in [0,x[\}\qquad V_x^{0-}=\inf\{n\geq 1; S_n\not\in [-x,0[\}$$
When applying the identity : $\cc B_q(x^-,iu) B_q(x,iu)- \cc C_q(x,iu)C_q(x^-,iu)=\phi(iu)+q$ to the compound 
Poisson process which of L\'evy exponent $1-\Q(e^{-iu S_1})$, on gets : For all $iu\in i\R$ and $s\in [0,1]$,

$$[1-\Q(e^{-iu S_{V_0^{x-}}}s^{V_0^{x-}};S_{V_0^{x-}}<0)][1-\Q(e^{-iu S_{V_x^{0-}}}s^{V_x^{0-}}; S_{V_x^{0-}}\geq 0)]$$
$$- \Q(e^{-iu S_{V_0^{x-}}}s^{V_0^{x-}}; S_{V_0^{x-}}\geq x)\Q(e^{-iu S_{V_x^{0-}}}s^{V_x^{0-}}; S_{V_x^{0-}}<-x)=1-s\Q(e^{-iu S_1})$$

\section{Examples} 
In this section, we treat two examples of L\'evy processes : stable processes which are not killed [$q=0$] and L\'evy processes without positive jumps. The "bilateral problem" has been essentially solved for these processes, by Rogozin \cite {R} in the first case and by Takacs \cite {T} in the second one. The reader can find also in  the recent work \cite{KK} 
 other  cases for which  this problem is solved.  Replacing in our context allows us to use the integral equations satisfied by our functions given in theorem \ref{ED}. This allows us to
  give some further properties, and  known results follow straightforwardly. 

\subsection{Stable processes } Let $\P$ be the law of a normalized stable process with index $\alpha$ and asymetry parameter $\rho$, (i;e. $\rho=\P(X_1>0)$).  Let us  $\gamma:=\alpha\rho$ and $\delta:=\alpha(1-\rho)$,  then $\gamma$ and $\delta$ belong to the interval $[0,1]$. The cases $\gamma=0$ and  $\delta=0$ corresponding to subordinators or opposite of subordinators are excluded in the sequel.
 The case $\gamma=1$ corresponds to stable processes without positive jumps.
The L\'evy exponent $\phi$ can be written as follows (we use the principal part of the  power functions) :
$$\phi(iu)=e^{ i{\pi\over 2}(\gamma-\delta)\epsilon (u)}|u|^{\alpha}=(-iu)^{\delta}(iu)^{\gamma}$$
We choose the spatial Wiener-Hopf factors 
$$\cc\psi_0(iu)=(-iu)^{\delta} \qquad \psi_0(iu)=(iu)^{\gamma} $$

\begin {theorem} For all $x\in ]0,+\i[$, one has 
$$A_0(x,\l)={1\over \Gamma(\gamma)} \int_0^x e^{-\l y} y^{\gamma-1} (1-{y\over x})^{\delta} dy\qquad (\l\in \C)$$
for  $\gamma\in ]0,1[$,
$$C_0(x,\l)={\Gamma(\delta+1)\over\Gamma(1-\gamma)\Gamma(\gamma)} \int_x^{+\i} e^{-\l y} ({y\over x}-1)^{-\gamma} y^{-\delta-1} dy\qquad (\Re(\l )\geq 0)$$
$$B_0(x,\l)= {\gamma \over \Gamma(1-\gamma)}\int _0^{+\i} (1-e^{-\l y} ) ({y\over x} +1)^{\delta}y^{-\gamma-1} dy  +  {\Gamma(\alpha)\over \Gamma(\delta)} x^{-\gamma}\qquad (\Re(\l)\geq 0)$$
and for  $\gamma=1$ (consequently  $\delta=\alpha-1$), 
$$C_0(x,\l)=\Gamma(\alpha).e^{-\l x} x^{-\delta}\qquad (\Re(\l)\geq 0)$$
$$B_0(x,\l)=\l + {\delta\over x} \qquad (\Re(\l)\geq 0)$$

The functions  $\cc A_0$, $\cc C_0$ and  $\cc B_0$ are obtained by duality.

\end{theorem} 

{\bf Sketch of proof} Stability property easily gives us that  the functions $H_0$ and $\cc H_0$ are of the form :
$$H_0(x)=k_+x^{\gamma}\qquad \cc H_0(x)=k_-x^{\delta}$$ 
where  $k_+$ and $k_-$ are positive constants.
Consequently, upon differentiating twice the integral equation of theorem \ref {ED}, one finds that the differential equation, in $x$,   satisfied by the three functions $A_0(x,\l)$, 
$C_0(x,\l)$ and  $\cc B_0(x,\l)$ is the hypergeometric confluent equation. This fact and the behavior of these functions when $x$ goes to $+\i$  allows us to compute them. We shall leave to the interested reader to check that one recovers in this way the results of Rogozin \cite{R} cited above.

\subsection{L\'evy processes without positive jumps}

In this section  $\P$ is the distribution of a L\'evy process without positive jumps. 
We refer to the book [B] chapter 7 for the results we remind here. The L\'evy exponent
is of the form 
$$\phi(iu)=\sigma^2u^2+aiu+\int_{]-\infty,0[}(1-e^{-iux}-iux\,1_{x>-1})\pi(dx)$$
This function can be extended to define  an analytic function on the half plane $\{\Re(\l)<0\}$ which is continous  when approaching the imaginary axis on the left. Still denote $\phi$ this function. Moreover, there exists a unique function $\Psi$ defined on $[0,+\i[$ such that $\phi(-\Psi(q))=-q$ for $q\geq 0$, and one can  take as  Wiener-Hopf factor $\psi_q$ the function : $$\psi_q(\l)=\l + \Psi(q)$$
This Wiener-Hopf factor corresponds to the choice of the local time $L^s(dt)=dS_t$. 
Notice that in case  $\P$ is the  law of a L\'evy process of the form $at-Y_t$ where $a$ is a positive real and $Y$ is a subordinator without drift, one has to take $L ^i(dt)={1\over a}\sum1_{I_t=X_t}\delta_t$ instead of $L ^i(dt)=\sum1_{I_t=X_t}\delta_t$ as we have stated at the beginning but this doesn't change previous results.

Denote $W_q$ the so called scale function : It is the unique increasing right continuous function defined for $q\in [0,+\i[$ and for $x\in[0,+\i[$ and satisfying the identity :
$$\int_0^{+\i} e^{-\l x} W_q(x) dx={1\over -\phi(-\l)-q}\qquad (\l>\Psi(q))$$
The associated Stieljes measure $W_q(dx)$ admits a right continous density $w_q$ on $]0,+\i[$ (we shall reprove this property in next proof):
$$W_q(x)=W_q(0)+\int _{]0,x]} w_q(y)dy$$
\begin{theorem} For  all $q \in [0,+\i[$ and $x\in ]0,+\i[$,  one has 
$$A_q(x,\l)=\int_0^x e^{-\l y}{W_q(x-y)\over W_q(x)}dy\qquad (\l\in\C)$$
$$\cc A_q(x,\l)=W_q(0)+\int_{-x}^0 e^{-\l y}[w_q(-y)- {w_q(x)\over W_q(x)}W_q(-y)] dy\qquad (\l\in\C)$$
The functions $C_q$ and $B_q$ can be extended to the whole complex plane and
$$ C_q(x,\l)= {e^{-\l x}\over W_q(x)} \qquad (\l\in\C)$$
$$B_q(x,\l)= \l + {w_q(x)\over W_q(x)}\qquad (\l\in\C)$$

$$\cc C_q(x,\l)= B_q(x,\l) -(\phi(\l)+q) \cc A_q(x,\l)\qquad (\Re(\l)\leq 0)$$
$$\cc B_q(x,\l)= (\phi(\l)+q) A_q(x,\l) + C_q(x,\l)\qquad (\Re(\l)\leq 0)$$
\end{theorem}
{\bf Sketch of proof} 
The absence of positive jumps gives us immediatly that 
$$C_q(x^-,\l)=e^{-\l x} N( S_{\xi_q}\geq x)\eqno{(5)}$$ $$B_q(x,\l)=\l +\cc N( \xi_q<\z \hbox{ or }  I_{\xi_q} <-x)\quad
 \hbox{and}\quad \lim_{x\to +\i}\cc N( \xi_q<\z \hbox{ or }  I_{\xi_q} <-x)=\Psi(q)$$
where $\xi_q$  denotes an independent exponential time after having enlarged the measure $N$ and $\cc N$ in order to contain such a variable. 

Derive  the integral equation  of theorem \ref {ED}  
$$C_q(x^-,\l)= \int_{[x,+\i[}e^{-\l y} B_q(y,\l) {H_q(dy)\over \cc H_q (y)}\eqno {(6)}$$ 
compare with the previous expression of $C_q$ and $B_q$ and deduce 
$${H_q(dx)\over \cc H_q (x)}=  N( S_{\xi_q}>x) dx \qquad \hbox{on }]0,+\i[\eqno{(7)}$$ 
 and the measure $N( S_{\xi_q}\in dx$) admits a  right continuous density $n_q$ on  $]0,+\i[$   and   
$$\cc N( \xi_q<\z \hbox{ or } I_{\xi_q} <-x)={n_q(x)\over N( S_{\xi_q}>x)} \qquad (x\in ]0,+\i[)$$
[In particular, $H_q(x)$, $N( S_{\xi_q}> x)$ and $C_q(x,\l)$ are continuous functions of $x$ on $]0,+\i[$]. 
 Consequently, the function ${n_q(x)\over N( S_{\xi_q}>x)}$ is nonincreasing and converges to $\Psi(q)$ when $x$ goes to  $+\i$.
 
 Remember  then the two following identities of theorem \ref {ED} :$$A_q(x^-,\l)= H_q(0)+\int_{]0,x]}e^{-\l y} \cc A_q(y,\l) {H_q(dy)\over \cc H_q (y)}\eqno {(8)}$$ 
$$B_q(x,\l)A_q(x^-,\l)+ C_q(x^-,\l)\cc A_q(x,\l)=1\eqno {(9)}$$
One first deduces from $(7)$ and $(8)$  that $A_q(x,\l)$ is a continuous function of $x$ (Notice that $A_q(x,0)=H_q(0)=0$ because $]0,+\i[$ is regular) and admits a right derivative. 
One then can compute the right derivative of the function ${A_q(x,\l)\over C_q(x,\l)}$ with the help of identities  $(5)$, $(6)$, $(7)$, $(8)$, $(9)$ and get :
$$\Bigl[{A_q(x,\l)\over C_q(x,\l)}\Bigr]'={B_q(x,\l)A_q(x^-,\l)+ C_q(x^-,\l)\cc A_q(x,\l)\over (e^{-\l x}N( S_{\xi_q}> x))^2} e^{-\l x}N( S_{\xi_q}> x)=e^{\l x} {1\over N( S_{\xi_q}> x)}$$
Integrate  this equation over $]0,x]$ and use equation (5) and  get : 
$$A_q(x,\l)= \int_0^x e^{-\l(x- y)}{N( S_{\xi_q}> x) \over N( S_{\xi_q}>y)} dy$$
Derive this equation  and use equation $(7)$ and $(8)$ and get  
$$\cc A_q(x,\l)= {1\over N( S_{\xi_q}>0)} +\int_0^x e^{\l y} ( {n_q(y)\over [N( S_{\xi_q}> y)]^2} - {1\over N( S_{\xi_q}>y)}  {n_q(x)\over N( S_{\xi_q}> x)})dy$$
For $\Re(\l)<0$, when  $x$  goes to $+\i$, $\cc A_q(x,\l)$ goes to ${1\over\cc \psi_q(\l)}= {\l+\Psi(q)\over \phi(\l)+q}$ and one gets :
$${\l+\Psi(q)\over \phi(\l)+q}={1\over N( S_{\xi_q}>0)}+\int_0^{+\i} e^{\l y} ( {n_q(y)\over [N( S_{\xi_q}>y)]^2} - {\Psi(q)\over N( S_{\xi_q}> y)}) dy$$
This equation allows us to identify   the  function ${1\over N( S_{\xi_q}>x)}$  to the scale function $W_q(x)$ and the statements of the theorem follow. \qed

\bibliographystyle{amsalpha}

\end{document}